\documentclass{gtart}

\def\ifplaintex{\expandafter\ifx\csname documentclass\endcsname\relax}

\ifplaintex 
\else
\fi

\expandafter\ifx\csname beginpicture\endcsname\relax
\expandafter\ifx\csname documentclass\endcsname\relax
\input pictex \else
\input prepictex \input pictex \input postpictex \fi\fi

\def\gt{{\mathsurround=0pt\it $\cal G\mskip-2mu$eometry \&\ 
$\cal T\!\!$opology}}        

\def\gtp{{\mathsurround=0pt\it $\cal G\mskip-2mu$eometry \&\ 
$\cal T\!\!$opology $\cal P\!$ublications}}  

\def\lognumber#1{\def\thelognumber{#1}}
\def\volumenumber#1{\def\thevolumenumber{#1}}
\def\papernumber#1{\def\thepapernumber{#1}}
\def\volumeyear#1{\def\thevolumeyear{#1}}

\def\pagenumbers#1#2{\def\startpage{#1}\def\finishpage{#2}}
\def\published#1{\def\publishdate{#1}}
\def\proposed#1{\def\theproposer{#1}}
\def\seconded#1{\def\theseconders{#1}}
\def\received#1{\def\receiveddate{#1}}
\def\revised#1{\def\reviseddate{#1}}
\def\accepted#1{\def\accepteddate{#1}}
\def\asciititle#1{\def\theasciititle{#1}}

\long\def\asciiabstract#1{\long\def\theasciiabstract{#1}}
\def\asciikeywords#1{\def\theasciikeywords{#1}}

\def\shorttitle#1{\def\theshorttitle{#1}}

\let\\\par\let\thelognumber\relax
\let\thevolumenumber\relax\let\thepapernumber\relax
\let\thevolumeyear\relax\let\thesamplenumber\relax\let\startpage\relax
\let\finishpage\relax\let\publishdate\relax\let\receiveddate\relax
\let\reviseddate\relax\let\accepteddate\relax\let\theasciititle\relax
\let\theasciiauthors\relax
\let\theasciiabstract\relax\let\theasciikeywords\relax
\let\theasciiemail\relax\let\theshortauthors\relax\let\theshorttitle\relax

\long\def\maketitlep{   

\count0=\startpage

\gt\hfill      
\beginpicture
\setcoordinatesystem units <0.33truein, 0.33truein> point at 2.2 0.9
\setplotsymbol ({$\cal G$})
\plotsymbolspacing=9truept
\circulararc 315 degrees from 0 1 center at 0 0
\setplotsymbol ({$\cal T$})
\circulararc 315 degrees from 1 -1 center at 1 0
\endpicture
\break
{\small\ifx\thesamplenumber\relax 
Volume \else Sample
\fi\thevolumenumber\ (\thevolumeyear)
\startpage--\finishpage\nl
Published: \publishdate}
\vglue 0.5truein plus 0.4fil minus 0.1truein

{\parskip=0pt\leftskip 0pt plus 1fil\def\\{\par\smallskip}{\ifplaintex\large
\else\Large\fi\bf\thetitle}\par\medskip}   

\vglue 0pt plus 0.1fil 

{\parskip=0pt\leftskip 0pt plus 1fil\def\\{\par}{\sc\theauthors}
\par\medskip}

\vglue 0pt plus 0.1fil 

{\small\parskip=0pt\let\newline\\
{\leftskip 0pt plus 1fil\def\\{\par}{\sl\theaddress}\par}
\expandafter\ifx\theemail\relax    
\relax\else\vglue 5pt plus 0.02fil minus 2pt\def\\{\stdspace{\rm 
and}\stdspace} 
\cl{Email:\stdspace\tt\theemail}\fi
\ifx\theurl\relax                  
\relax\else\vglue 5pt plus 0.02fil minus 2pt\def\\{\stdspace{\rm 
and}\stdspace}
\cl{URL:\stdspace\tt\theurl}\fi\par}

\vglue 7pt plus 0.3fil minus 3pt

{\bf Abstract}
\vglue 5pt plus 0.1fil minus 2pt

\theabstract

\vglue 7pt plus 0.3fil minus 3pt

{\bf AMS Classification numbers}\quad Primary:\quad \theprimaryclass

Secondary:\quad \thesecondaryclass

\vglue 5pt plus 0.3fil minus 2pt

{\bf Keywords}\quad \thekeywords

\vglue 10pt plus 0.5fil minus 5pt

{\small  Proposed: \theproposer\hfill Received: \receiveddate\nl
Seconded: \theseconders\hfill 
\ifx\reviseddate\relax                         
Accepted: \accepteddate                        
\else
Revised: \reviseddate                          
\fi}
\eject
}       

\let\maketitlepage\maketitlep
\let\maketitle\maketitlepage

\font\phead=cmsl9 scaled 950
\font=cmsl9 scaled 1050
\font\pnum=cmbx10 scaled 913
\font\lnum=cmbx10 
\font\pfoot=cmsl9 scaled 950
\font=cmsl9 scaled 1050
\ifplaintex
\headline{\vbox to 0pt{\vskip -4.5mm\line{\small\phead\ifnum
\count0=\startpage ISSN 1364-0380 (on line)
1465-3060 (printed) \hfill {\pnum\folio}\else\ifodd\count0\def\\{ }%
\ifx\theshorttitle\relax\thetitle\else\theshorttitle\fi\hfill{\pnum\folio}
\else\def\\{ and }{\pnum\folio}\hfill\ifx\theshortauthors\relax\theauthors
\else\theshortauthors\fi\fi\fi}\vss}}
\footline{\vbox to 0pt{\vglue 0mm\line{\small\pfoot\ifnum\count0=\startpage
\copyright\ \gtp\hfill\else
\gt, Volume \thevolumenumber\ (\thevolumeyear)\hfill\fi}\vss
}}
\else
\makeatletter
\def\@oddhead{{\small\lhead\ifnum\count0=\startpage ISSN 1364-0380 (on line)
1465-3060 (printed) \hfill {\lnum\number\count0}\else\ifodd\count0
\def\\{ }\ifx\theshorttitle\relax \thetitle \else\theshorttitle\fi\hfill
{\lnum\number\count0}\else\def\\{ and }{\lnum\number\count0}
\hfill\ifx\theshortauthors\relax 
\theauthors\else\theshortauthors\fi\fi\fi}}\def\@evenhead{\@oddhead}
\def\@oddfoot{\small\lfoot\ifnum\count0=\startpage\copyright\ \gtp\hfill\else
\gt, Volume \thevolumenumber\ (\thevolumeyear)\hfill\fi}
\def\@evenfoot{\@oddfoot}
\makeatother
\fi

\newwrite\gtoutfile
\long\gdef\makeheadfile{  
{\def\\{, }\def\s{ }
\immediate\openout\gtoutfile head.xxx
\immediate\write\gtoutfile{To: math@arxiv.org}
\immediate\write\gtoutfile{Subject: put or rep NNNNN:pppp}
\immediate\write\gtoutfile{--text follows this line--}
\immediate\write\gtoutfile{Proxy-for: \ifx\theasciiauthors\relax
\theauthors\else\theasciiauthors\fi\s<\ifx\theasciiemail\relax\theemail\else\theasciiemail\fi>}
\immediate\write\gtoutfile{\noexpand\\}
\immediate\write\gtoutfile{Authors: \ifx\theasciiauthors\relax
\theauthors\else\theasciiauthors\fi}
{\def\\{ }\immediate\write\gtoutfile{Title: \ifx\theasciititle\relax
\thetitle\else\theasciititle\fi}}
\immediate\write\gtoutfile{Subj-class: GT or SG or MG etc}
\immediate\write\gtoutfile{MSC-class: \theprimaryclass\ifx\thesecondaryclass\relax\else, \thesecondaryclass\fi}
\immediate\write\gtoutfile{Journal-ref: Geom. Topol. \thevolumenumber
(\thevolumeyear) \startpage-\finishpage}
\immediate\write\gtoutfile{Comments: Published by Geometry and Topology at}
\immediate\write\gtoutfile{\s\s http://www.maths.warwick.ac.uk/gt/GTVol\thevolumenumber/paper\thepapernumber.abs.html}
\immediate\write\gtoutfile{\noexpand\\}
\immediate\write\gtoutfile{}
\ifx\theasciiabstract\relax
\immediate\write\gtoutfile{\theabstract}\else
\immediate\write\gtoutfile{\theasciiabstract}\fi
\immediate\write\gtoutfile{}
\immediate\write\gtoutfile{\noexpand\\}
\immediate\write\gtoutfile{}
\immediate\closeout\gtoutfile}}  

\def\maketitlepage{\maketitlep\makeheadfile}
\let\maketitle\maketitlepage

\lognumber{203}
\volumenumber{7}\papernumber{13}
\volumeyear{2003}
\pagenumbers{443}{486}
\received{19 August 2001}
\revised{7 June 2003}
\accepted{10 July 2003}
\published{18 July 2003\nl  Republished with corrections: 21 August 2003}
\proposed{Walter Neumann}
\seconded{Benson Farb, Martin Bridson}

\usepackage{amsmath,amscd,amssymb,epsf,subfigure}
\usepackage[mathscr]{eucal}
\newtheorem{thm}{Theorem}
\newtheorem*{thm*}{Theorem}
\newtheorem{lemma}[thm]{Lemma}

\newtheorem{proposition}[thm]{Proposition}

\newtheorem*{corollary*}{Corollary}
\newtheorem*{claim*}{Claim}

\numberwithin{equation}{subsection}
\numberwithin{thm}{subsection}
\begin{document}
\newcommand{\R}{{\mathbb R}}
\newcommand{\C}{{\mathbb C}}
\newcommand{\Z}{{\mathbb Z}}
\newcommand{\B}{{\mathbb B}}
\newcommand{\Id}{{\mathbb I}}
\renewcommand{\P}{{\mathbb P}}
\newcommand{\cpo}{{\C\P^1}}
\newcommand{\X}{{\mathfrak X}}
\newcommand{\E}{{\mathfrak E}}
\newcommand{\Ss}{{\mathfrak S}}
\newcommand{\M}{{\pi_0(\o{Homeo}(M))}}
\newcommand{\Mb}{{\pi_0(\o{Homeo}(M,\partial M))}}
\renewcommand{\o}{\operatorname}
\newcommand{\Hom}{\o{Hom}}
\newcommand{\GLtz}{{\GL{2,\Z}}}
\newcommand{\PGLtz}{{\o{PGL}(2,\Z)}}
\newcommand{\Tei}{{\mathfrak T}}
\newcommand{\Teich}{\Tei_M}
\newcommand{\kt}{{\kappa^{-1}(t)}}
\newcommand{\rr}{{\mathcal{R}}}
\newcommand{\ktR}{{\kappa^{-1}(t)\cap\R^3}}
\newcommand{\ktwo}{{\kappa^{-1}(2)}}
\newcommand{\ktwoR}{{\ktwo\cap\R^3}}
\newcommand{\intt}{{\o{int}}}
\newcommand{\Hh}{{\mathfrak H}}
\newcommand{\Cc}{{\mathfrak C}}
\newcommand{\Hm}[1]{\Hom (\pi ,#1)}
\newcommand{\Hmm}[1]{\Hom (\pi ,#1)/#1}
\newcommand{\Hmmm}[1]{\Hom (\pi ,#1)/\hspace{-3pt}/#1}
\newcommand{\hmg}{\Hom (\pi ,G)}
\renewcommand{\k}{\mathbf{k}}
\newcommand{\Ker}{\o{Ker}}
\newcommand{\Aut}{\o{Aut}}
\newcommand{\Ad}{\o{Ad}}
\newcommand{\Inn}{\o{Inn}}
\newcommand{\Out}{\o{Out}}
\newcommand{\Ou}{\Out({\pi})}
\newcommand{\hg}{\Hmmm{G}}
\newcommand{\tr}{\o{tr}}
\newcommand{\dd}[1]{\frac{\partial}{\partial{#1}}}
\newcommand{\SL}[1]{{\o{SL}}({#1})}
\newcommand{\SLt}{{\SL{2}}}
\newcommand{\GL}[1]{{\o{GL}}({#1})}
\newcommand{\slt}{{\SL{2,\C}}}
\newcommand{\slr}{{\SL{2,\R}}}
\newcommand{\sothc}{{\o{SO}{(3,\C)}}}
\newcommand{\soth}{{\o{SO}(3)}}
\newcommand{\soto}{{\o{SO}(2,1)}}
\newcommand{\soot}{{\o{SO}(1,2)}}
\newcommand{\sovb}{{\o{SO}(\C^3,\B)}}
\newcommand{\sovbr}{{\o{SO}(\R^3,\B)}}
\newcommand{\sot}{{\o{SO}{(2)}}}
\newcommand{\sooo}{{\o{SO}{(1,1)}}}
\newcommand{\Ht}{\mathbf{H}^2}
\newcommand{\rpt}{{\R\P}^2}
\newcommand{\sut}{{\o{SU}}(2) }
\newcommand{\uo}{{\o{U}}(1)}
\newcommand{\cho}{{\mathbf{H}}^1_{\C}}
\newcommand{\qq}{{\mathbf Q}}
\newcommand{\qqq}{{\mathcal Q}}
\newcommand{\qd}{{\ddddot{\mathcal Q}}}
\newcommand{\tu}{{\tilde u}}
\newcommand{\tx}{{\tilde x}}
\newcommand{\ty}{{\tilde y}}
\newcommand{\tz}{{\tilde z}}
\newcommand{\tg}{{\tilde\gamma}}
\newcommand{\Ak}{{\Aut(\C^3,\kappa)}}
\renewcommand{\ll}{{\mathcal L}}
\def\S{Section }

\title{The modular group action on real $SL(2)$--characters\\of a
one-holed torus}
\asciititle{The modular group action on real SL(2)-characters of a
one-holed torus}
\shorttitle{Action of the modular group}

\author{William M Goldman}
\address{Mathematics Department,
University of Maryland\\College Park, MD  20742 USA}
\email{wmg@math.umd.edu}



\asciikeywords{Surface, fundamental group, character variety,
representation variety, mapping class group, ergodic action, proper
action, hyperbolic structure with cone singularity, Fricke space,
Teichmueller space}

\begin{abstract}
The group $\Gamma$ of automorphisms of the 
polynomial 
$$\kappa(x,y,z) = x^2 + y^2 + z^2 - x y z -2 $$
is isomorphic to 
$$\PGLtz\ltimes(\Z/2\oplus\Z/2).$$
For $t\in\R$, the $\Gamma$-action  on $\ktR$ displays rich and
varied dynamics.  The action of $\Gamma$ preserves a Poisson structure
defining a $\Gamma$--invariant area form on each $\ktR$. For
$t<2$, the action of $\Gamma$ is properly discontinuous on the four
contractible components of $\ktR$ and ergodic on the compact component
(which is empty if $t<-2$). The contractible components correspond to
Teichm\"uller spaces of (possibly singular) hyperbolic structures on a
torus $\bar{M}$. 
For $t=2$, the level set $\ktR$ consists of characters of reducible 
representations and comprises two ergodic components
corresponding to actions of $\GLtz$ on $(\R/\Z)^2$ and $\R^2$
respectively.  For $2 < t \le 18$, the action of $\Gamma$ on $\ktR$ is
ergodic.  Corresponding to the Fricke space of a three-holed 
sphere is a $\Gamma$--invariant open subset $\Omega\subset\R^3$ 
whose components 
are permuted freely by a subgroup of index $6$ in $\Gamma$.
The level set $\ktR$
intersects $\Omega$ if and only if $t>18$, in which case 
the $\Gamma$--action on the complement $(\ktR)-\Omega$ is ergodic.
\end{abstract}
\asciiabstract{The group Gamma of automorphisms of the polynomial
kappa(x,y,z) = x^2 + y^2 + z^2 - xyz -2 is isomorphic to PGL(2,Z)
semi-direct product with (Z/2+Z/2).  For t in R, Gamma-action on ktR =
kappa^{-1}(t) intersect R displays rich and varied dynamics.  The
action of Gamma preserves a Poisson structure defining a
Gamma-invariant area form on each ktR. For t < 2, the action of Gamma
is properly discontinuous on the four contractible components of ktR
and ergodic on the compact component (which is empty if t < -2). The
contractible components correspond to Teichmueller spaces of (possibly
singular) hyperbolic structures on a torus M-bar.  For t = 2, the
level set ktR consists of characters of reducible representations and
comprises two ergodic components corresponding to actions of GL(2,Z)
on (R/Z)^2 and R^2 respectively.  For 2 < t <= 18, the action of Gamma
on ktR is ergodic.  Corresponding to the Fricke space of a three-holed
sphere is a Gamma-invariant open subset Omega subset R^3 whose
components are permuted freely by a subgroup of index 6 in Gamma.  The
level set ktR intersects Omega if and only if t > 18, in which case
the Gamma-action on the complement ktR - Omega is ergodic.}

{\small\maketitlepage}

\section*{Introduction}

Let $M$ be a compact oriented surface of genus one with one boundary
component, a {\em one-holed torus.\/}  Its fundamental group $\pi$ is free
of rank two.  Its mapping class group $\M$ is isomorphic to the outer
automorphism group $\Ou$ of $\pi$ and acts on the space of equivalence
classes of representations $\pi\longrightarrow\SLt$. We investigate the
dynamics of this action on the set of real points on this moduli space.

Corresponding to the boundary of $M$ is an element $K\in\pi$ which is the
commutator of free generators $X,Y\in\pi$. By a theorem of 
Fricke~\cite{Fricke,FrickeKlein}, the moduli space of 
equivalence classes of $\slt$--representations
naturally identifies with affine 3--space $\C^3$, via the quotient map
\begin{align*}
\Hm{\slt} & \longrightarrow  \C^3 \\
\rho & \longmapsto 
\bmatrix x \\ y  \\ z \endbmatrix  :=
\bmatrix \tr (\rho (X))\\
\tr (\rho (Y))\\ \tr (\rho (XY))\endbmatrix.
\end{align*}
In terms of these coordinates, the trace $\tr\rho(K)$ equals:
\begin{equation*}
\kappa(x,y,z) := x^2 + y^2 + z^2 - xyz -2
\end{equation*}
which is preserved under the action of $\Ou$.  The action of $\Ou$ on
$\C^3$ is commensurable with the action of the group $\Gamma$ of
polynomial automorphisms of $\C^3$ which preserve $\kappa$
(Horowitz~\cite{Horowitz}). (Compare also Magnus~\cite{Magnus}.)

\begin{thm*} Let 
$\kappa(x,y,z) = x^2 + y^2 + z^2 - xyz -2$ and let $t\in\R$.
\begin{itemize}
\item For $t<-2$, the group $\Gamma$ acts properly on
$\ktR$;
\item For $-2 \le  t < 2$, there is a compact connected component $C_t$
of $\kappa^{-1}(t)\,\cap$ $\R^3$  upon which $\Gamma$ acts ergodically; 
$\Gamma$ acts properly on the complement $\big(\ktR\big)- C_t$;
\item For $t = 2$, the action of $\Gamma$ is ergodic on the compact
subset $\kappa^{-1}(2)\cap [-2,2]^3$ and the action is ergodic on
the complement $\kappa^{-1}(2)- [-2,2]^3$;
\item For $18 \ge t > 2$, the group $\Gamma$ acts ergodically on 
$\ktR$;
\item For $t> 18$, the group $\Gamma$ 
acts properly and freely on an open subset 
$\Omega_t\subset\ktR$, permuting its components. The $\Gamma$--action on the 
complement of $\Omega_t$ is ergodic.
\end{itemize}
\end{thm*}
The proof uses the interplay between representations of the fundamental
group and hyperbolic structures on $M$. The dynamics breaks up into
two strikingly different types: representations corresponding to
hyperbolic structures comprise contractible connected components of
the the level sets $\kt$, whereas representations which map a simple
nonperipheral essential loop to an elliptic element, comprise open
subsets of $\kt$ upon which $\Gamma$ is ergodic. Thus nontrivial
dynamics accompanies nontrivial topology of the moduli spaces.

In his doctoral thesis~\cite{Stantchev}, 
G.\ Stantchev considers characters corresponding to representations
into $\o{PGL}(2,\R)$ (that is, actions preserving $\Ht$ but not
preserving orientation on $\Ht$). For $t< -14$ (respectively 
$t > 6$)  characters of discrete
embeddings representing hyperbolic structures on 2--holed projective planes
(respectively 1--holed Klein bottles)
give wandering domains in the corresponding level set. For $ -14 \le t < 2$
the $\Gamma$--action on the corresponding level set is ergodic, and for
$t < -14$, the action is ergodic on the complement of the wandering domains
corresponding to Fricke spaces of 2--holed projective planes.

\begin{figure}[ht!] 
\cl{
\subfigure[Level set $\kappa=-2.1$]{\epsfxsize=2.3in \epsfbox{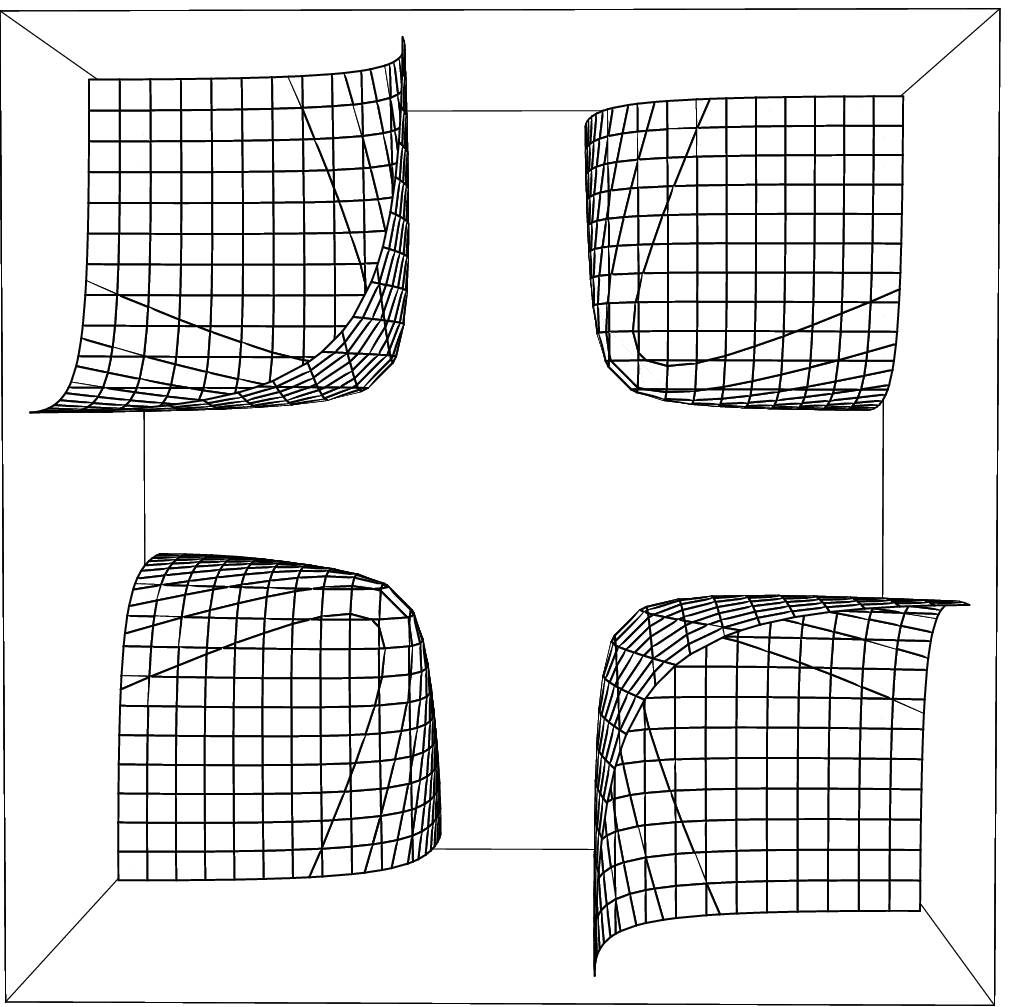}}
\qquad
\subfigure[Level set $\kappa=1.9$]{\epsfxsize=2.3in \epsfbox{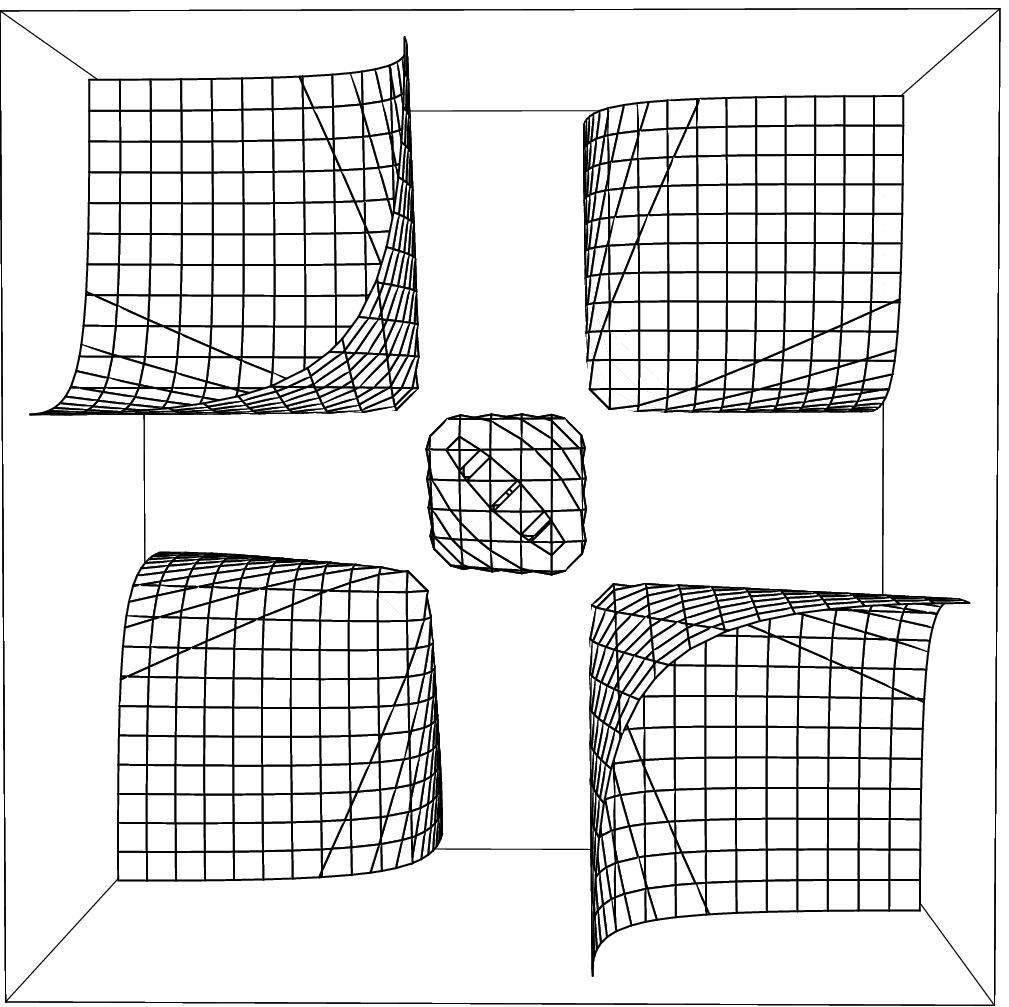}}
	}\vspace{-5mm}
\nocolon\caption{}
\end{figure}

\begin{figure}[ht!] 
\cl{
\subfigure[Level set $\kappa=1.9$]{\epsfxsize=2.3in \epsfbox{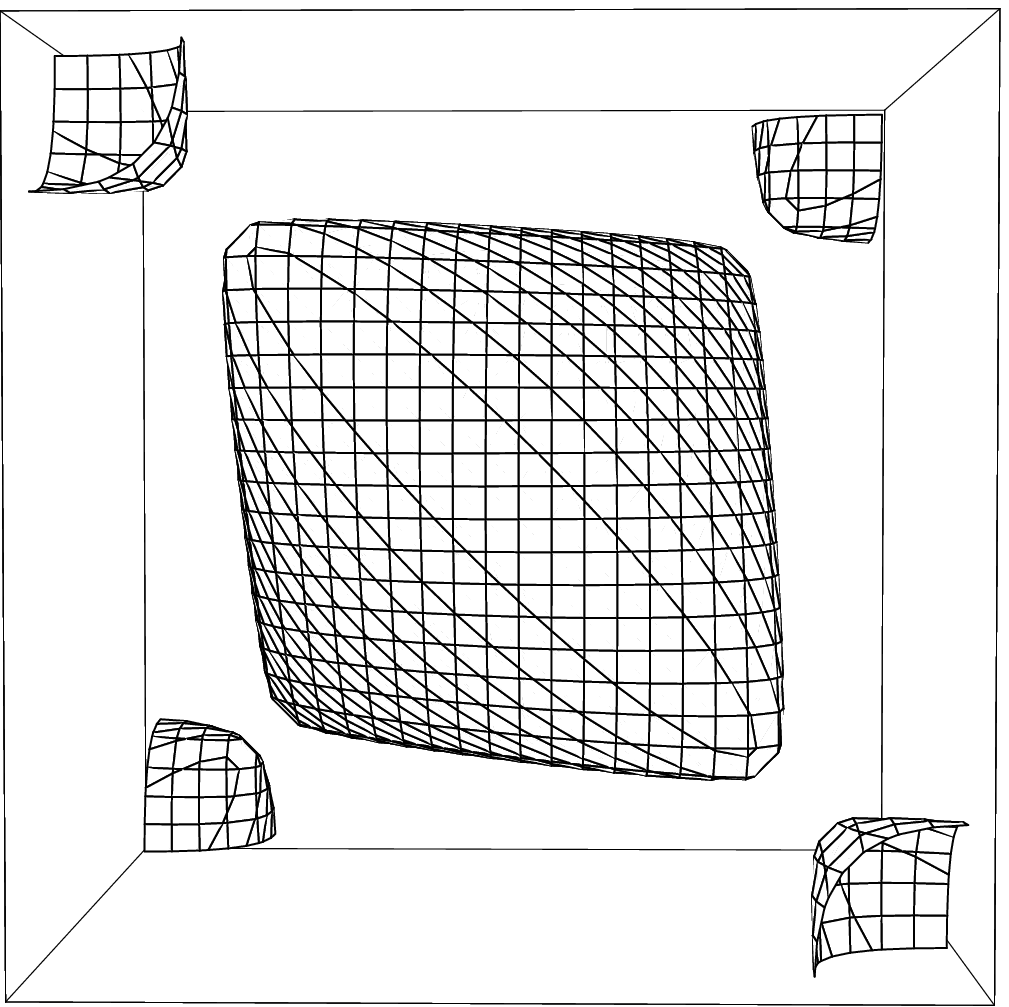}}
\qquad
\subfigure[Level set $\kappa=2.1$]{\epsfxsize=2.3in \epsfbox{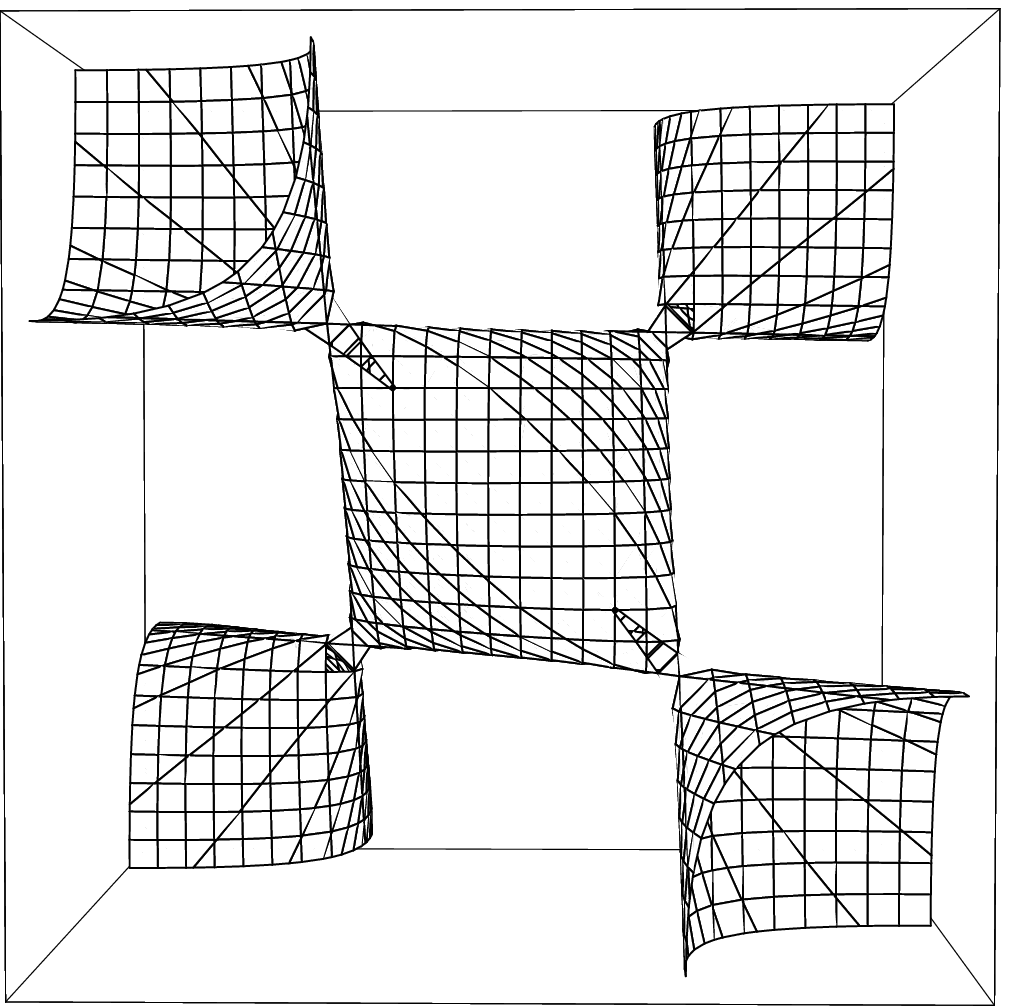}}
	}\vspace{-5mm}
\nocolon\caption{}
\end{figure}

\begin{figure}[ht!] 
\cl{
\subfigure[Level set $\kappa=2.1$]{\epsfxsize=2.3in \epsfbox{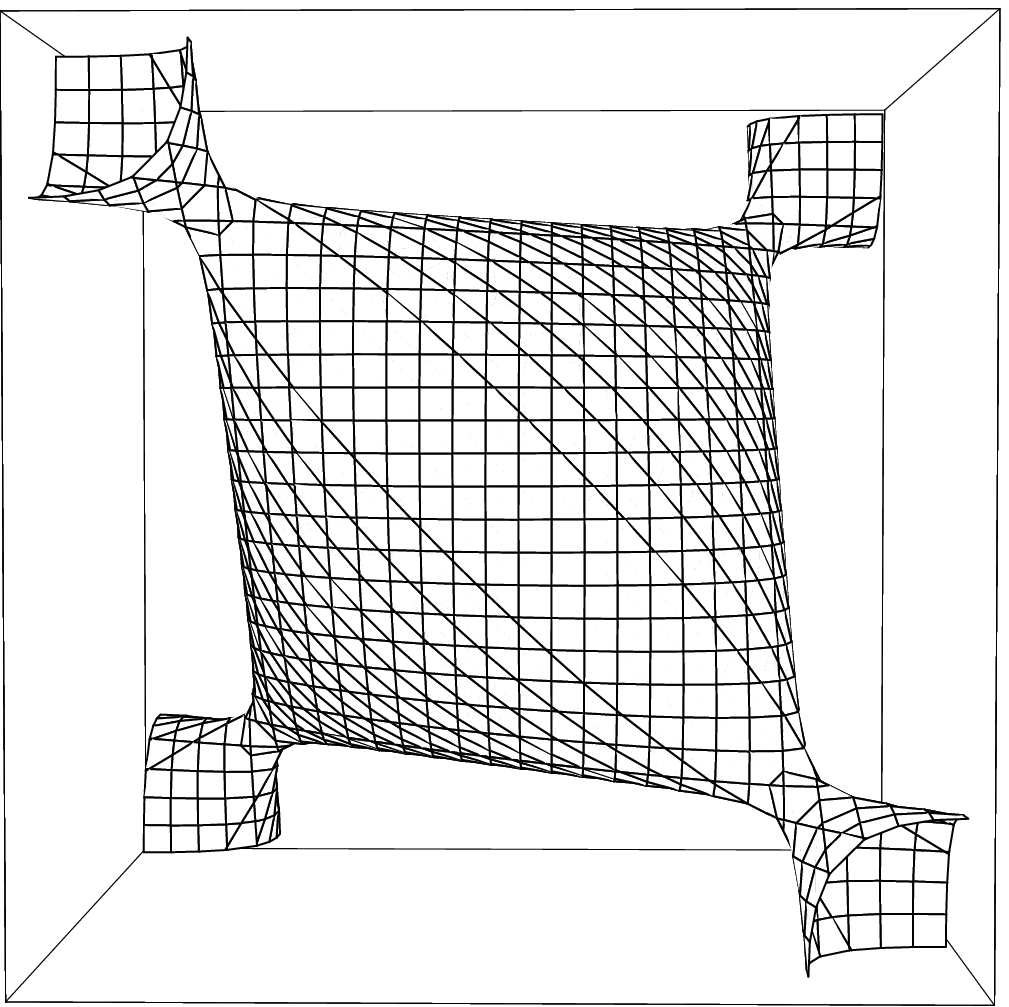}}
\qquad
\subfigure[Level set $\kappa=4$]{\epsfxsize=2.3in \epsfbox{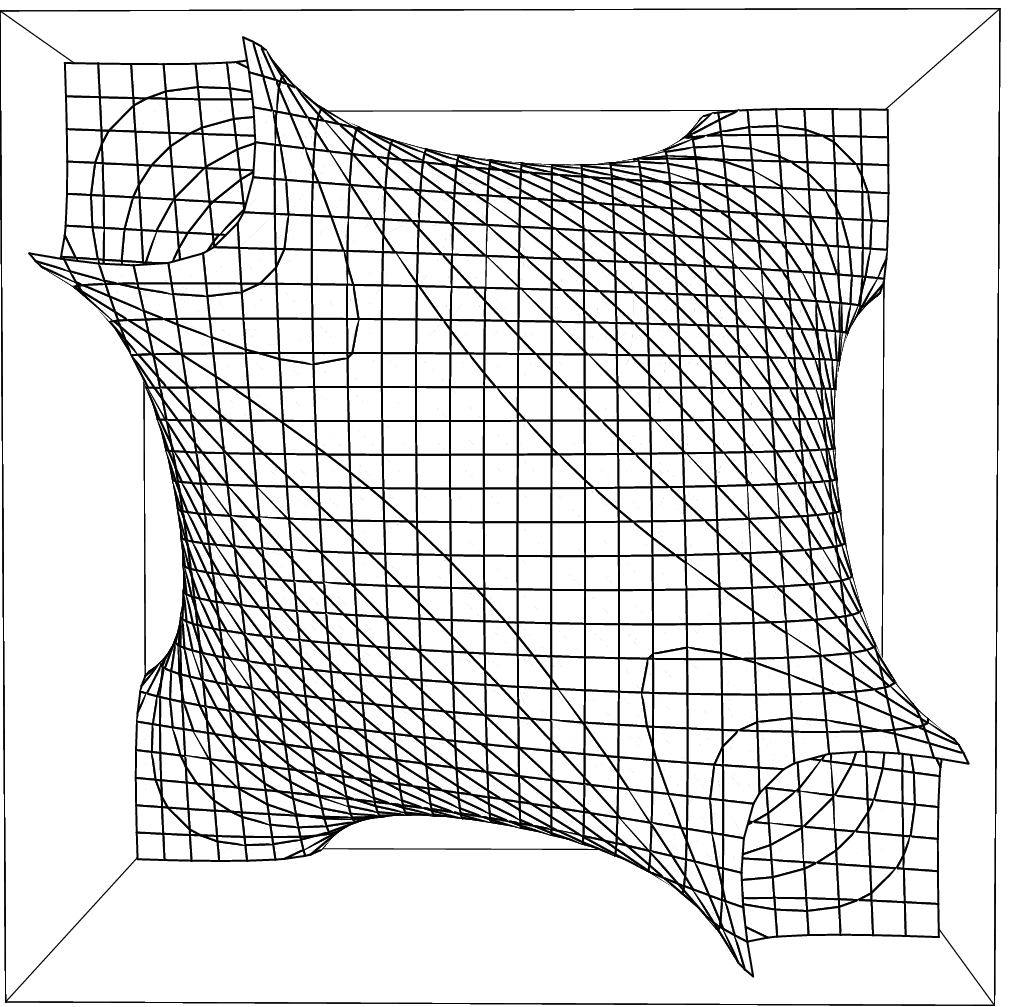}}
	}\vspace{-5mm}
\nocolon\caption{}
\end{figure}

\begin{figure}[ht!] 
\cl{
\subfigure[Level set $\kappa=4$]{\epsfxsize=2.3in \epsfbox{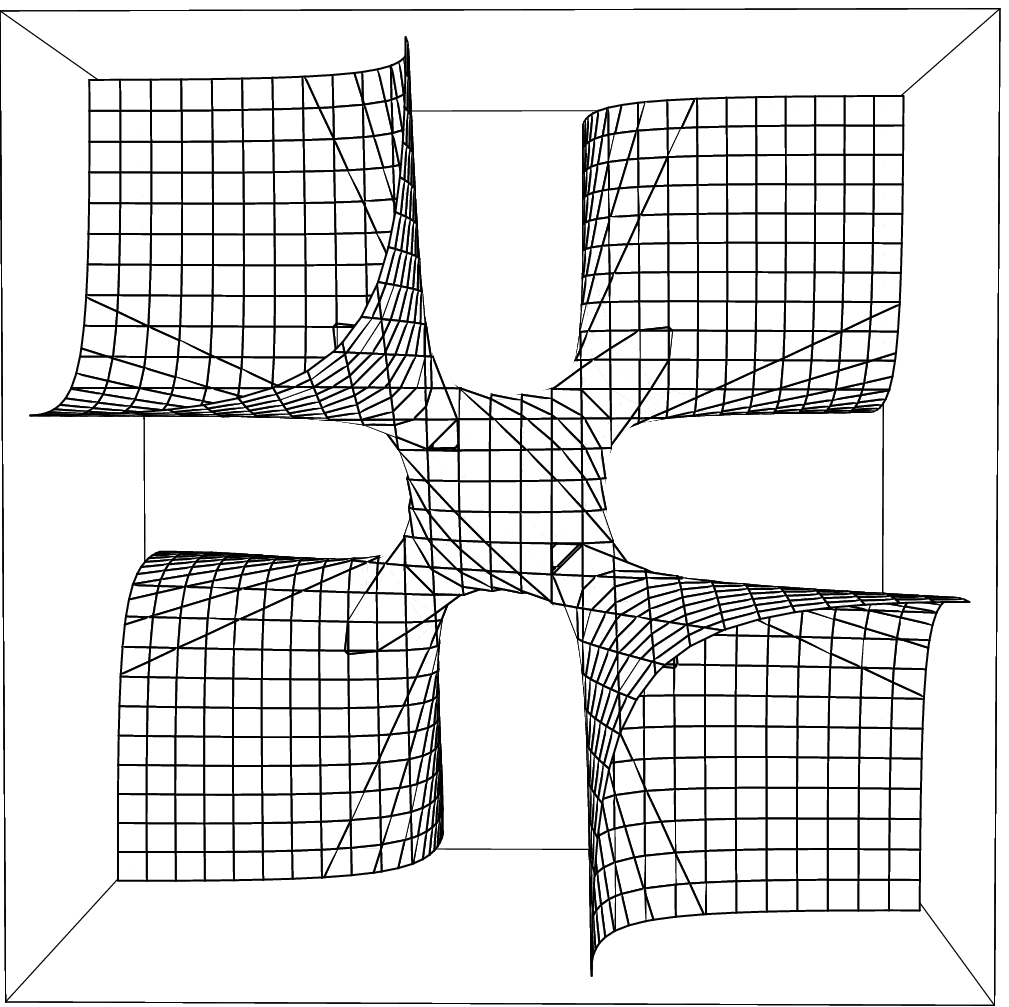}}
\qquad
\subfigure[Level set $\kappa=19$]{\epsfxsize=2.3in \epsfbox{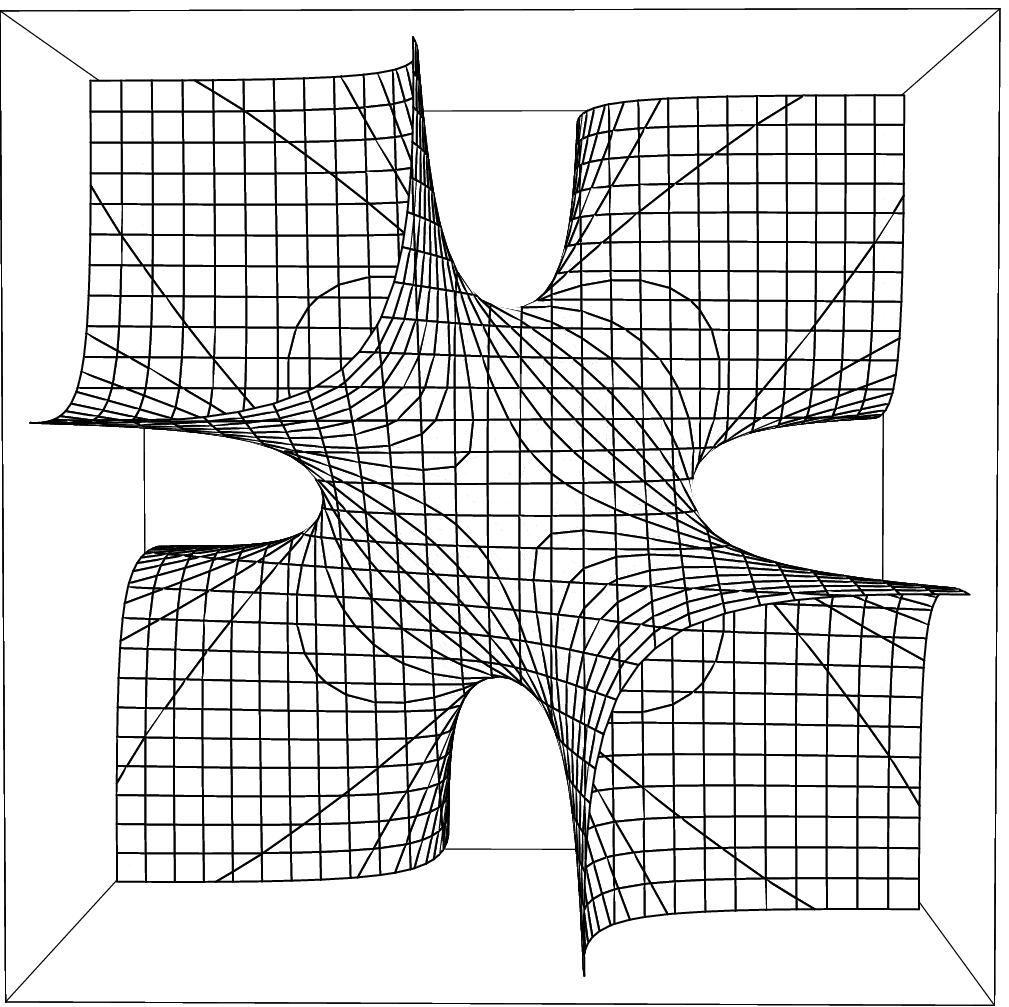}}
	}\vspace{-5mm}
\nocolon\caption{}
\end{figure}

{\bf Notation and terminology}\qua
We work in Poincar\'e's model of the hyperbolic plane $\Ht$ as the
upper half-plane. We denote the identity map (identity matrix) by $\Id$. 
Commutators are denoted: 
\begin{equation*}
[A,B] = ABA^{-1}B^{-1}  
\end{equation*}
and inner automorphisms are denoted:
\begin{equation*}
\iota_g\co  x \longmapsto g x g^{-1}. 
\end{equation*}
For any group $\pi$ we denote the group of all automorphisms 
by $\Aut(\pi)$ and the normal subgroup of inner automorphisms 
by $\Inn(\pi)$. We denote the quotient group $\Aut(\pi)/\Inn(\pi)$ by 
$\Out(\pi)$.
A compact surface $S$ with $n$ boundary components will be called 
{\em $n$--holed.\/} Thus a {\em one-holed torus\/} is the complement
of an open disc inside a torus and a {\em three-holed sphere\/} 
(sometimes called a ``pair-of-pants'') is the complement of three disjoint
discs inside a sphere.

\medskip
{\bf Acknowledgements}\qua  I am grateful to Matt Bainbridge, Robert
Benedetto, Lawrence Ein, Cathy Jones, Misha Kapovich, Bernhard Leeb,
John Millson, John Parker, Dan Rudolph, Peter Shalen,
Richard Schwartz, Adam Sikora, and Scott Wolpert for helpful
conversations during the course of this work.
I especially want to thank Walter Neumann and George Stantchev for
carefully reading the manuscript and suggesting major improvements.

I gratefully acknowledge partial support from National Science
Foundation grant DMS--9803518, DMS--0103889 and a Semester Research
Award from the General Research Board of the University of Maryland in
Fall 1998.

\section{The modular group and the moduli space}
In this section we define the modular group and the moduli space for the
one-holed torus. The modular group is isomorphic to 
\begin{equation*}
\M\cong\Ou\cong\GLtz 
\end{equation*}
and the moduli space is affine space $\C^3$. We explain how the invariant
function $\kappa$ originates.
We describe various elements of the modular group and how they act on the
moduli space.
\subsection{The mapping class group}
Let $M$ denote a compact connected orientable surface of genus one
with one boundary component. Since attaching a disc to $M$ yields a
torus, we refer to $M$ as a {\em one-holed torus.\/} 
The {\em mapping class group\/} of $M$ is the group $\M$ 
of isotopy classes of homeomorphisms of $M$. 
We investigate the action of $\M$ on the moduli space of
flat $\SLt$--connections on $M$.
\subsubsection{Relation to $\pi_1(M)$}
Choose a basepoint $x_0\in M$ and let
$\pi := \pi_1(M;x_0)$.
Any homeomorphism of $M$ is isotopic to one which fixes $x_0$
and hence defines an automorphism of $\pi$.
Two such isotopic homeomorphisms determine
automorphisms of $\pi$ differing by an inner automorphism, 
producing a well-defined homomorphism
\begin{equation}\label{eq:Nielsen}
N\co \M \longrightarrow \Ou \cong \Aut(\pi)/\Inn(\pi).
\end{equation}
If $M$ is a closed surface, then
Dehn (unpublished) and Nielsen~\cite{Nielsen1927} proved
that $N$ is an isomorphism. (See Stillwell~\cite{Stillwell} 
for a proof of the Dehn--Nielsen theorem.) 
When $\partial M\neq\emptyset$, then each component $\partial_i M$
determines a conjugacy class $\mathcal{C}_i$ of elements of $\pi_1(M)$.
The image of $N$ consists of elements of $\Ou$ represented by
automorphisms which preserve each $\mathcal{C}_i$.

Let $M$ be a one-holed torus. Its fundamental group admits
a geometric redundant presentation
\begin{equation}\label{eq:pi1torus}
\pi = \langle X,Y,K \mid [X,Y] = K \rangle
\end{equation}
where $K$ corresponds to the generator of $\pi_1(\partial M)$.
Of course, $\pi$ is freely generated by $X,Y$.

\subsubsection{Nielsen's theorem}
The following remarkable property of $M$ is due to
Jakob Nielsen~\cite{Nielsen} and does not generalize to other
hyperbolic surfaces with boundary.  
For a proof see  Magnus--Karrass--Solitar~\cite{MKS}, Theorem~3.9
or Lyndon--Schupp~\cite{LyndonSchupp}, Proposition 5.1.

\begin{proposition}\label{prop:Nielsen}
Any automorphism of the rank two group
\begin{equation*}
\pi=\langle X,Y,K\mid K= [X,Y]\rangle  
\end{equation*}
takes $K$ to a conjugate of either $K$ itself or its inverse
$K^{-1}$.
\end{proposition}
An equivalent geometric formulation is:
\begin{proposition}\label{prop:NielsenGeometric}
Every homotopy-equivalence $M\to M$ is homotopic to a
homeomorphism of $M$.  
\end{proposition}
Thus the homomorphism \eqref{eq:Nielsen}
defines an isomorphism 
\begin{equation*}
N\co \M\longrightarrow\Ou. 
\end{equation*}

\subsection{The structure of the modular group}\label{sec:structure}

We say that an automorphism of $\pi$ which takes $K$ to either $K$ or
$K^{-1}$ is {\em normalized.\/} The normalized automorphisms form a
subgroup $\Aut(\pi,K)$ of $\Aut(\pi)$. 

Let $\phi$ be an automorphism of $\pi$. Proposition~\ref{prop:Nielsen} implies
that $\phi(K)$ is conjugate to $K^{\pm 1}$. Thus 
there exists an inner automorphism $\iota_g$ such that
\begin{equation*}
\phi(K) = gK^{\epsilon}g^{-1},
\end{equation*}
that is 
\begin{equation*}
\iota_g^{-1}\circ\phi\in\Aut(\pi,K)
\end{equation*}
where $\epsilon = \pm 1$. (In fact, $\epsilon = \det(h([\phi]))$,
where $h$ is the homomorphism defined by \eqref{eq:homology} below.)
Since the centralizer of $K$ in $\pi$ equals the cyclic group $\langle
K\rangle$, the automorphism $\phi$ determines the coset of $g$ modulo
$\langle K\rangle$ uniquely.

We obtain a short exact sequence
\begin{equation*}
1 \longrightarrow \langle \iota_K\rangle
  \hookrightarrow \Aut(\pi,K) 
  \longrightarrow \Ou
  \longrightarrow 1
\end{equation*}
The action on the homology $H_1(M;\Z)\cong \Z^2$ defines a homomorphism
\begin{equation}\label{eq:homology}
h\co \Ou \longrightarrow \GLtz.
\end{equation}
By Nielsen~\cite{Nielsen}, $h$ is an isomorphism.  
(Surjectivity follows by realizing an 
element of $\GLtz$ as a {\em linear\/} homeomorphism of the torus
$\R^2/\Z^2$.  See Lyndon--Schupp~\cite{LyndonSchupp}, Proposition 4.5 or
Magnus--Karrass--Solitar~\cite{MKS}, \S3.5, Corollary N4.)  We obtain
an isomorphism 
\begin{equation*}
\mu\co \GLtz \longrightarrow \Aut(\pi,K)/\langle\iota_K\rangle.
\end{equation*}
Restriction of the composition $\mu\circ h$ to $\Aut(\pi,K)$
equals the quotient homomorphism
\begin{equation*}
\Aut(\pi,K)\longrightarrow\Aut(\pi,K)/\langle\iota_K\rangle.
\end{equation*}

\section{Structure of the character variety}

\subsection{Trace functions}

The relevant moduli space is the {\em character variety,\/}
the categorical quotient of $\hmg$ by the $G$--action
by inner automorphisms, where $G=\slt$. Since
$\pi$ is freely generated by two elements $X$ and $Y$, the set 
$\hmg$ of homomorphisms $\pi\longrightarrow G$ identifies with the set of
pairs $(\xi,\eta)\in G\times G$, via the mapping
\begin{align*}
\hmg & \longrightarrow G\times G \\
\rho & \longmapsto (\rho(X),\rho(Y)).
\end{align*}
This mapping is equivariant respecting the action of $G$ on $\hmg$ by
\begin{equation*}
g\co  \rho \longmapsto \iota_g\circ\rho 
\end{equation*}
and the action of $G$ on $G\times G$ by
\begin{equation*}
g\co (\xi,\eta) \longmapsto (g\xi g^{-1},g\eta g^{-1}).
\end{equation*}
The moduli space $\hg$ consists of equivalence classes of elements of
$\hmg\cong G\times G$ where the equivalence class of a homomorphism $\rho$ 
is defined as the closure of the $G$--orbit $G\rho$. Then $\hg$ is the
categorical quotient in the sense that its coordinate ring identifies
with the ring of $G$--invariant regular functions on $\hmg$.
For a single element $g\in G$, the conjugacy class $\Inn(G)(g)$ is determined
by the trace $t = \tr(g)$ if $t\neq\pm 2$. That is, 
\begin{equation} \label{eq:singleconjugacyclass}
\Inn(G)(g) = \tr^{-1}(t). 
\end{equation}
For $t=\pm 2$, then
\begin{equation}\label{eq:singleconjugacyclass2}
\tr^{-1}(\pm 2) = \{\pm \Id\}\quad  \bigcup \quad \Inn(G)\cdot\left( 
\pm \bmatrix 1 & 1 \\ 0 & 1\endbmatrix\right).
\end{equation}
By Fricke~\cite{Fricke} and Fricke--Klein~\cite{FrickeKlein}, 
the traces of the generators $X,Y,XY$ param\-et\-rize
$\hg$ as the affine space $\C^3$. 
As $X$ and $Y$ freely generate $\pi$, we may identify:
\begin{align*}
\hmg  & \longleftrightarrow \quad G\times G \\
\rho \quad & \longleftrightarrow (\rho(X), \rho(Y))
\end{align*}
The character mapping
\begin{align*}
\chi\co  \hg & \longrightarrow \C^3 \\
[\rho]  & \longmapsto 
\bmatrix x(\rho) \\ y(\rho)  \\ z(\rho) \endbmatrix  = 
\bmatrix \tr (\rho (X))\\
\tr (\rho (Y))\\ \tr (\rho (XY))\endbmatrix
\end{align*}
is an isomorphism. (Compare the discussion in Goldman~\cite{TopComps},
4.1, \cite{Erg}, Sections 4--5 and \cite{Expo}.)

For example, given $(x,y,z)\in\C^3$, the representation
$\rho$ defined by 
\begin{equation}\label{eq:explicitrep}
\rho(X) = \bmatrix x &  -1  \\ 1 & 0 \endbmatrix,\quad
\rho(Y) = \bmatrix 0 &  \zeta^{-1}  \\ -\zeta & y \endbmatrix
\end{equation}
satisfies $\chi(\rho)=(x,y,z)$ where $\zeta\in\C$ is chosen so that
\begin{equation*}
\zeta  + \zeta^{-1}=  z. 
\end{equation*}
Conversely, if $(x,y,z)\in\C^3$ then $\chi^{-1}(x,y,z)$ consists of a
single $G$--orbit if and only if $\kappa(x,y,z)\neq 2$ where $\kappa$
is defined below in \eqref{eq:defkappa}.  (This is also the condition that
$\rho$ is an {\em irreducible representation\/} of $\pi$. Compare
Lubotzky--Magid~\cite{LubotzkyMagid},
Brumfiel--Hilden~\cite{BrumfielHilden},
Culler--Shalen~\cite{CullerShalen}.)

For any word $w(X,Y)$, the function 
\begin{align*}
\hmg & \longrightarrow \C \\
[\rho] & \longmapsto \tr(\rho(w(X,Y))
\end{align*}
is $G$--invariant. Hence there exists a polynomial $f_w(x,y,z)\in\C[x,y,z]$
such that
\begin{equation*}
\tr(\rho(w(X,Y)) = f_w(x(\rho),y(\rho),z(\rho)).
\end{equation*}
A particularly important example occurs for $w(X,Y)=[X,Y] = K$,
in which case we denote $f_w(x,y,z)$ by $\kappa(x,y,z)$. 
By an elementary calculation (see, for example~\cite{Expo}),
\begin{equation}\label{eq:defkappa}
\tr(\rho(K)) = \kappa(x,y,z) := x^2 + y^2 + z^2 - xyz -2.
\end{equation}
The level set $\kt$ consists of equivalence classes of representations
$\rho\co \pi\to G$ where $\rho(K)$ is constrained to lie in 
$\tr^{-1}(t)$. (Compare \eqref{eq:singleconjugacyclass} and
\eqref{eq:singleconjugacyclass2} above.)

\subsection{Automorphisms}
Let $G = \slt$.
The group $\Aut(\pi)$ acts on the character variety\break $\hg$ by:
\begin{equation*}
\phi_*([\rho]) = [\rho\circ\phi^{-1}]
\end{equation*}
and since
\begin{equation*}
\rho\circ\iota_\gamma = \iota_{\rho(\gamma)}\circ\rho
\end{equation*}
the subgroup $\Inn(\pi)$ acts trivially. Thus $\Ou$ acts on
\begin{equation*}
\hg\cong\C^3 
\end{equation*}
and since an automorphism $\phi$ of $\pi$ is determined by 
\begin{equation*}
(\phi(X),\phi(Y))= (w_1(X,Y),w_2(X,Y)),  
\end{equation*}
the action of $\phi$ on
$\C^3$ is given by the three polynomials $f_{w_1},f_{w_2},f_{w_3}$:
\begin{equation*}
\bmatrix x \\ y \\ z \endbmatrix \longmapsto
\bmatrix f_{w_1}(x,y,z) \\ f_{w_2}(x,y,z) \\ f_{w_3}(x,y,z)
\endbmatrix.
\end{equation*}
where $w_3(X,Y) = w_1(X,Y)w_2(X,Y)$ is the product word.
Hence $\Ou$ acts on $\C^3$ by {\em polynomial automorphisms.\/}
Nielsen's theorem (Proposition~\ref{prop:Nielsen}) 
implies that any such automorphism preserves $\kappa\co \C^3\longrightarrow\C$,
that is
\begin{equation*}
\kappa\left(f_{w_1}(x,y,z),f_{w_2}(x,y,z),f_{w_1w_2}(x,y,z)\right) 
= \kappa(x,y,z).
\end{equation*}

\subsubsection{Sign-change automorphisms}\label{sec:signchanges}
Some automorphisms of the character variety are not induced
by automorphisms of $\pi$. 
Namely, the homomorphisms of $\pi$ into the center
$\{\pm\Id\}\subset G$ form a group acting
on $\hmg$ by pointwise multiplication. 
Let $\zeta\in\Hom(\pi,\{\pm\Id\})$ and $\rho\in\hmg$. 
Then 
\begin{equation}\label{eq:twisted}
\zeta\cdot\rho\co  \gamma \longmapsto \zeta(\gamma)\rho (\gamma) 
\end{equation}
is a homomorphism. 
This defines an action
\begin{equation*}
\Hom(\pi,\{\pm\Id\}) \times \hmg \longrightarrow \hmg.
\end{equation*}
Furthermore, since $\{\pm\Id\}$ 
is central in $G$ and $K\in\pi$ is a commutator
\begin{equation}\label{eq:twistkom}
(\zeta\cdot\rho)(K) = \rho(K). 
\end{equation}
Since $\pi$ is free of rank two, 
\begin{equation*}
\Hom(\pi,\{\pm\Id\})  \cong \Z/2 \times \Z/2. 
\end{equation*}
The three nontrivial elements 
$(0,1),(1,0),(1,1)$  
of $\Z/2 \times \Z/2$  act on representations
by $\sigma_1,\sigma_2,\sigma_3$ respectively:
\begin{equation*}
\sigma_1\cdot\rho\co  
\begin{cases}
X & \longmapsto \rho(X) \\
Y & \longmapsto -\rho(Y) \end{cases}
\end{equation*}
\begin{equation*}
\sigma_2\cdot\rho\co  
\begin{cases}
X & \longmapsto -\rho(X) \\
Y & \longmapsto \rho(Y) \end{cases}
\end{equation*}
\begin{equation*}
\sigma_3\cdot\rho\co  
\begin{cases}
X & \longmapsto -\rho(X) \\
Y & \longmapsto -\rho(Y) \end{cases}
\end{equation*}
The corresponding action on characters is:
\begin{align*}
(\sigma_1)_*:\bmatrix x \\ y \\ z \endbmatrix
& \longmapsto \bmatrix x \\ -y \\ -z \endbmatrix \\
(\sigma_2)_*:\bmatrix x \\ y \\ z \endbmatrix
& \longmapsto \bmatrix -x \\ y \\ -z \endbmatrix \\
(\sigma_3)_*:\bmatrix x \\ y \\ z \endbmatrix
& \longmapsto \bmatrix -x \\ -y \\ z \endbmatrix. 
\end{align*}
We call this group the {\em group of sign-changes\/} and denote it by $\Sigma$.

Evidently $\Sigma$ preserves $\kappa(x,y,z)$. (Alternatively apply
\eqref{eq:twistkom}.)

\subsubsection{Permutations}\label{sec:permutations}

Since $\kappa(x,y,z)$ is symmetric in $x,y,z$ the full symmetric group
$\Ss_3$ also acts on $\C^3$ preserving $\kappa$. Unlike $\Sigma$, 
elements of $\Ss_3$ are induced by automorphisms of $\pi$. 
The group of all {\em linear\/} 
automorphisms of $(\C^3,\kappa)$ is generated by $\Sigma$ and $\Ss_3$ and
forms a semidirect product $\Sigma \rtimes \Ss_3$.

$\Ss_3$ is actually a quotient of $\Gamma$.
The {\em projective line $\P^1(\Z/2)$ over $\Z/2$\/} has
three elements, and every permutation of this set is realized
by a projective transformation. Since $\Z/2$ has only one nonzero element,
the projective automorphism group of $\P^1(\Z/2)$ equals $\GL{2,\Z/2}$
and $\GL{2,\Z/2} \cong \Ss_3$.
The action
of $\GLtz$ on $(\Z/2)^2$ defines a homomorphism
\begin{equation*}
\GLtz  \longrightarrow \GL{2,\Z/2} \cong \Ss_3
\end{equation*}
whose kernel $\GLtz_{(2)}$ is generated by the involutions
\begin{equation}\label{eq:gens2congsubgp}
\bmatrix 1 & 0 \\ 0 & -1 \endbmatrix, 
\bmatrix 1 & -2 \\ 0 & -1 \endbmatrix, 
\bmatrix 1 & 0 \\ 2 & -1 \endbmatrix 
\end{equation}
and $-\Id\in\GLtz_{(2)}$. The sequence
\begin{equation*}
\PGLtz_{(2)} :=  \GLtz_{(2)}/\{\pm\Id\} \longrightarrow
\PGLtz \longrightarrow \Ss_3
\end{equation*}
is exact.
The kernel $\Gamma_{(2)}$ of the composition
\begin{equation*}
\Gamma \longrightarrow \PGLtz \longrightarrow  \Ss_3,
\end{equation*}
equals the  semidirect product
\begin{equation*}
\PGLtz_{(2)}\ltimes(\Z/2\oplus\Z/2)
\end{equation*}
and $\Gamma$ is an extension
\begin{equation*}
\Gamma_{(2)} \longrightarrow \Gamma \longrightarrow \Ss_3.
\end{equation*}
(See Goldman--Neumann~\cite{GoldmanNeumann} for more extensive discussion of 
$\Gamma$ and its action on the the set of {\em complex points\/} of the
character variety.)

\subsubsection{Other automorphisms}\label{sec:other}
For later use, here are several specific elements of $\Gamma$.
See also the appendix to this paper.

The {\em elliptic involution \/} is the automorphism
\begin{align*}
X & \longmapsto  X^{-1} \\
Y & \longmapsto  Y^{-1} 
\end{align*}
corresponding to $-\Id\in\GLtz$ and acts trivially on characters.

The {\em Dehn twist\/} about $X$ is the automorphism
$\tau_X\in\Aut(\pi)$
\begin{align*}
X & \longmapsto  X \notag \\
Y & \longmapsto  Y X \notag 
\end{align*}
inducing the automorphism of characters
\begin{equation*} 
\bmatrix x \\ y \\ z \endbmatrix \longmapsto
\bmatrix x \\ xy - z \\ y \endbmatrix
\end{equation*}
and corresponds to
\begin{equation*}
\bmatrix 1 & 1 \\ 0 & 1 \endbmatrix\in\GLtz. 
\end{equation*}
The {\em Dehn twist\/} about $Y$ is the automorphism
$\tau_Y\in\Aut(\pi)$
\begin{align*}
X & \longmapsto  X Y \notag \\
Y & \longmapsto  Y  \notag 
\end{align*}
inducing the automorphism of characters
\begin{equation*} 
\bmatrix x \\ y \\ z \endbmatrix \longmapsto
\bmatrix x y - z \\ y \\ x \endbmatrix.
\end{equation*}
and corresponds to 
\begin{equation*}
\bmatrix 1 & 0 \\ 1 & 1 \endbmatrix\in\GLtz.  
\end{equation*}
The {\em quadratic reflection $Q_z$\/} is the automorphism
\begin{align*}
X & \longmapsto  X  \notag \\
Y & \longmapsto  Y^{-1} \notag 
\end{align*}
inducing the automorphism of characters
\begin{equation*} 
\bmatrix x \\ y \\ z \endbmatrix \longmapsto
\bmatrix x \\ y \\  x y - z \endbmatrix.
\end{equation*}
and corresponds to 
\begin{equation*}
\bmatrix 1 & 0 \\ 0 & -1 \endbmatrix\in\GLtz.  
\end{equation*}

\subsection{Reducible characters}\label{sec:reducible}
Any representation $\rho$ having character in $\ktwo$ is reducible,
as can be checked by direct calculation of commutators in $\slt$.
Namely, let $\xi,\eta\in\slt$. We may write
\begin{equation*}
\eta = \bmatrix a  & b \\ c  & d \endbmatrix 
\end{equation*}
where $ad-bc = 1$. By applying an inner automorphism, we may assume
that $\xi$ is in Jordan canonical form. If $\xi$ is diagonal,
\begin{equation*}
\xi = \bmatrix \lambda  & 0 \\ 0 & \lambda^{-1} \endbmatrix
\end{equation*}
then
\begin{equation*}
\tr [\xi,\eta] = 2 + bc (\lambda - \lambda^{-1})^2 
\end{equation*}
implies that if $\tr[\xi,\eta]=2$, then either $\xi=\pm\Id$ or
$bc=0$ (so $\xi$ is upper-triangular or lower-triangular).
Otherwise 
\begin{equation*}
\xi = \pm \bmatrix 1 & s \\ 0 & 1 \endbmatrix,\quad 
\end{equation*}
(where $s\neq 0$), in which case
\begin{equation}\label{eq:trcomunip}
\tr [\xi,\eta] = 2 + s^2 c^2 \ge 2
\end{equation}
so $\tr[\xi,\eta]=2$ implies that $c=0$ and $\eta$ is upper-triangular.
Thus if $\tr(\rho(K))=2$, then 
$\rho$ is conjugate to an upper-triangular representation.

We may replace $\rho$ by its {\em semisimplification,\/} that is
the upper-triangular matrices by the corresponding diagonal matrices
(their semisimple parts) to obtain a representation by diagonal
matrices having the same character:
\begin{align}\label{eq:uppertri}
\rho(X)  & = \bmatrix \xi & * \\ 0 & \xi^{-1} \endbmatrix \\
\rho(Y)  & = \bmatrix \eta & * \\ 0 & \eta^{-1} \endbmatrix \notag \\
\rho(XY)  & = \bmatrix \zeta & * \\ 0 & \zeta^{-1} \endbmatrix. \notag 
\end{align}
where $\xi\eta=\zeta$.
Thus
\begin{align}\label{eq:extension}
x & = \xi + \xi^{-1}, \\  
y & = \eta + \eta^{-1},\notag \\   
z & = \zeta + \zeta^{-1}.\notag
\end{align}
\eqref{eq:extension} corresponds to the following factorization 
\eqref{eq:factorization} of 
\begin{equation*}
\kappa(x,y,z) - 2 = x^2 + y^2 +z^2 - xyz -4.    
\end{equation*}
Under the embedding 
\begin{equation*}
\C[x,y,z] \longrightarrow \C[\xi,\xi^{-1},\eta,\eta^{-1},\zeta,\zeta^{-1}]
\end{equation*}
defined by \eqref{eq:extension}, the polynomial $\kappa(x,y,z)-2$ factors:
\begin{align}\label{eq:factorization}
\kappa( & \xi + \xi^{-1}, \eta + \eta^{-1},\zeta + \zeta^{-1} ) - 2  
=  \\
& \zeta^{-2} 
(1 - \xi\eta\zeta) 
(1 - \xi^{-1}\eta\zeta)
(1 - \xi\eta^{-1}\zeta)  
(1 - \xi\eta\zeta^{-1}). \notag
\end{align}
Given $(x,y,z)$ with $\kappa(x,y,z)=2$, the triple $(\xi,\eta,\zeta)$
is only defined up to an action of  $\Z/2\oplus\Z/2$. Namely the automorphisms
\begin{equation*}
\begin{cases}
X & \longmapsto X^{-1} \\
Y & \longmapsto Y \end{cases}\qquad
\begin{cases}
X & \longmapsto X \\
Y & \longmapsto Y^{-1} \end{cases}\qquad
\begin{cases}
X & \longmapsto X^{-1} \\
Y & \longmapsto Y^{-1} \end{cases}
\end{equation*}
act on $(\xi,\eta,\zeta)$ by:
\begin{equation*}
\begin{cases}
\xi & \longmapsto \xi^{-1} \\
\eta & \longmapsto \eta  \\
\zeta & \longmapsto \zeta^{-1}\end{cases}\qquad\qua
\begin{cases}
\xi & \longmapsto \xi \\
\eta & \longmapsto \eta^{-1}  \\
\zeta & \longmapsto \zeta^{-1}\end{cases}\qquad\qua
\begin{cases}
\xi & \longmapsto \xi^{-1} \\
\eta & \longmapsto \eta^{-1}  \\
\zeta & \longmapsto \zeta\end{cases}
\end{equation*}
respectively.

\subsection{The Poisson structure}
As in Goldman~\cite{Erg},\S 5.3, the automorphisms preserving $\kappa$
are unimodular and therefore preserve the exterior bivector field as
well
\begin{align*}
\Xi  & =   \frac1{2\pi^2} d\kappa \cdot \dd{x}\wedge\dd{y}\wedge\dd{z} \\
& = \frac1{2\pi^2} 		
\Bigg( \ 
   (2x - yz) \dd{y}\wedge\dd{z}  \\ 
& \qquad \qquad +  (2y - zx) \dd{z}\wedge\dd{x} \\
& \qquad \quad \qquad+  (2z - xy) \dd{x}\wedge\dd{y} \quad \Bigg) 
\end{align*}
which restricts to (the dual of) an area form on each level set $\kt$
which is invariant under $\Ou$. We shall always consider this invariant
measure on $\ktR$. 
In the case when $\kt$ has a rational parametrization by an affine plane, 
that is, when $t=2$, the above bivector field is the image of the constant
bivector field under the parametrization:
\begin{align*}
\C\times \C & \longrightarrow \, \ktwo \\
(\xi,\eta) & \longmapsto 
\bmatrix e^{\xi} +  e^{-\xi} \\ e^{\eta} +  e^{-\eta} \\
e^{\xi + \eta} +  e^{-(\xi+\eta)} \endbmatrix \\
\dd{\xi}\wedge\dd{\eta}  & \longmapsto  \,2\pi^2 
\Xi
\end{align*}
The bivector field $\Xi$ defines a {\em Poisson structure\/} on the
moduli space $\C^3$ for which $\kappa$ defines a {\em Casimir function.\/}

(The coefficient $2\pi^2$ occurs to normalize the area of the 
compact component $C_K$ of reducible $\sut$--characters to $1$.)
Here we study the action of $\Ou$ on the set $\ktR$ of $\R$--points of 
$\kt$ where $t\in\R$, with respect to this invariant measure.

\subsection{The orthogonal representation}\label{sec:orthogonal}
To describe the structure of the character variety more completely,
we use a 3--dimensional orthogonal representation associated to
a character $(x,y,z)\in\C^3$. (Compare \S 4.2 of Goldman~\cite{TopComps}
and Brumfiel--Hilden~\cite{BrumfielHilden}.) This will be used to
identify real characters as characters of representations into the real
forms $\slr$ and $\sut$ of $\slt$.

Consider the complex vector space $\C^3$ with the standard basis 
$\{e_1,e_2,e_3\}$ and bilinear form
defined by the symmetric matrix
\begin{equation*}
\B = \bmatrix 2 & z & y \\ z & 2 & x \\ y & x & 2 \endbmatrix.  
\end{equation*}
Since 
\begin{equation*}
\det(\B) = -2 \big(\kappa(x,y,z) - 2\big) 
\end{equation*}
(where $\kappa(x,y,z)$ is defined in \eqref{eq:defkappa}), 
the symmetric bilinear form is nondegenerate if and only if 
$\kappa(x,y,z)\neq 2$. 

Assume $\kappa(x,y,z)\neq 2$ so that $\B$ is nondegenerate. 
Let $\sovb\cong\sothc$ denote the group of unimodular linear transformations 
of $\C^3$ orthogonal with respect to $\B$. The local isomorphism
\begin{equation*}
\Phi\co  \slt \longrightarrow \sovb
\end{equation*}
is a surjective double covering, equivalent to the adjoint
representation (or the representation on the second symmetric power of $\C^2$),
and is unique up to composition with automorphisms of $\slt$ and $\sothc$.

For any vector $v\in \C^3$ such that $\B(v,v)\neq 0$,  the reflection
\begin{equation}\label{eq:reflection}
R_v\co  x \longmapsto x - 2\;\frac{\B(v,x)}{\B(v,v)} v
\end{equation}
is a $\B$--orthogonal involution.

Let $\Z/2\star\Z/2\star\Z/2$ be the free product of cyclic
groups $\langle R_i\rangle$, where $i=1,2,3$, and $R_i^2=\Id$. The homomorphism
\begin{align*}
\pi & \longrightarrow \Z/2\star\Z/2\star\Z/2  \\
X &\longmapsto \quad R_2 R_3\\
Y &\longmapsto \quad R_3 R_1
\end{align*}
embeds $\pi$ as an index two subgroup of $\Z/2\star\Z/2\star\Z/2$.

\begin{lemma}
Let $\rho\in\hmg$ satisfy $\kappa(x,y,z)\neq 2$. The restriction to $\pi$
of the representation $\hat\rho\co \Z/2\star\Z/2\star\Z/2
\longrightarrow\sovb$ defined by $R_i \longmapsto R_{e_i}$
is conjugate to the composition
$\Phi\circ\rho\co \pi\longrightarrow\sovb $.
\end{lemma}
See \S 4.2 of Goldman~\cite{TopComps},
or Brumfiel--Hilden~\cite{BrumfielHilden} for a proof.

\subsection{The set of $\R$--points}\label{sec:realpoints}

A real character is the character of a representation into a {\em real
form\/} of $G$. (Proposition III.1.1 of Morgan--Shalen~\cite{MorganShalen}).
The above orthogonal representation gives an explicit form of this result.

Suppose $(x,y,z)\in\R^3$ and $\kappa(x,y,z)\neq 2$. Then the
restriction of $\B$ to $\R^3$ is a nondegenerate symmetric $\R$--valued
bilinear form (also denoted $\B$), which is either indefinite or
(positive) definite.  Let $\sovbr$ denote the group of unimodular
linear transformations of $\R^3$ orthogonal with respect to $\B$.

There are two conjugacy classes of real forms of $G$,
compact and noncompact.  Every compact real form of $G$ is conjugate to
$\sut$ and every noncompact real form of $G$ is conjugate to $\slr$.
Specifically, if $\B$ is positive definite, then $\sovbr$ is conjugate
to $\soth$ and $\Phi^{-1}(\sovbr)$ is conjugate to $\sut$.
If $\B$ is indefinite, then  $\sovbr$ is conjugate to either
$\soto$ or $\soot$ (depending on whether $\kappa(x,y,z)>2$ or
$\kappa(x,y,z)<2$ respectively) and 
$\Phi^{-1}(\sovbr)$ is conjugate to $\slr$.
Thus every real character $(x,y,z)\in\R^3$ is the character of a representation
into either $\sut$ or $\slr$. 

When $x,y,z\in\R$, the restriction of $\B$ to $\R^3$ is definite
if and only if $-2\le x,y,z\le 2$ and $\kappa(x,y,z)<2$.

The coordinates of characters of representations
$\pi_1(M)\longrightarrow\sut$ satisfy:
\begin{equation*}
-2\le x,y,z \le 2,\quad \quad x^2 + y^2 + z^2 - xyz -2 	\le 2	
\end{equation*}
and, for fixed $t\in[-2,2]$, comprise a component of
$\ktR$.
The other four connected components of $\ktR$ consist of characters of
$\slr$--representations. These four components are 
freely permuted by $\Sigma$.  For $t<-2$, all four components of the
relative character variety $\ktR$ consist of characters of
$\slr$--representations and are freely permuted by $\Sigma$.  
For $t>2$, the relative character variety
$\ktR$ is connected and consists of characters of
$\slr$--representations.

The two critical values $\pm 2$ of $\tr\co G\longrightarrow\C$ deserve
special attention. When $t=-2$, a representation $\rho\in\hmg$ with
$\tr\rho(K)=-2$ is a regular point of the mapping
\begin{align*}
E_K\co \hmg & \longrightarrow G \\
\rho \qquad & \longmapsto \rho(K).
\end{align*}
Such a representation is a regular point of the composition
\begin{equation*}
\tr\circ E_K\co   \hmg  \longrightarrow \C
\end{equation*}
unless $\rho(K)= -I$. In the latter case (when $[\rho(X),\rho(Y)]=-I$),
the representation is conjugate to the {\em quaternion representation
in $\sut$\/}
\begin{equation}\label{eq:quaternionrep}
\rho(X) = \bmatrix i & 0 \\ 0 & -i\endbmatrix, \quad
\rho(Y) = \bmatrix 0 & -1 \\ 1 & 0\endbmatrix. 
\end{equation}
In both cases, $G$ acts locally freely on the subset
$(\tr\circ E_K)^{-1}(-2)$  of $\hmg$ with quotient $\kappa^{-1}(-2)$.

Now consider the case $t=2$.
As in \S\ref{sec:reducible}, the $\R$--points of $\ktwo$ correspond to 
reducible representations, and in fact are characters of 
representations with values 
in the Cartan subgroups (maximal tori) 
of the real forms $\sut$ and $\slr$.
Every Cartan subgroup of $\sut$ is conjugate to $\uo$, and
every Cartan subgroup of $\slr$ is conjugate to either
$\sot$ or $\sooo$.
Characters of reducible $\sut$--representations form the compact set
\begin{equation*}
C_K = \ktwo \cap [-2,2]^3 
\end{equation*}
which identifies with the quotient of the 2--torus $\uo\times\uo$
by $\{\pm\Id\}$ under the extension \eqref{eq:extension}. 

Characters of reducible $\slr$--representations comprise four components
(related by $C_i = \sigma_i C_0$ for $i=1,2,3$):
\begin{align*}
C_0 & = 
\ktwo  \cap \bigg( [2,\infty)\times[2,\infty)\times[2,\infty)\bigg), \\ 
C_1 & =
\ktwo \cap \bigg( [2,\infty)\times(-\infty,-2]\times(-\infty,-2]\bigg), \\ 
C_2 & =
\ktwo \cap \bigg( (-\infty,-2] \times[2,\infty)\times(-\infty,-2]\bigg),\\ 
C_3 & =
\ktwo \cap \bigg( (-\infty,-2] \times (-\infty,-2]\times[2,\infty)\bigg),
\end{align*}
each of which identifies with the quotient of $\R_+\times\R_+$
by $\{\pm\Id\}$ under \eqref{eq:extension}. These four components are
freely permuted by $\Sigma$.
There is a compact component of equivalence classes of $\slr$--representations
which are irreducible but not {\em absolutely irreducible\/} --- that is,
although $\R^2$ is an irreducible $\pi$--module, 
its complexification is reducible.
In that case the representation is conjugate to a representation in
$\sot$. This component agrees with the component 
$C_K$ consisting of characters of reducible $\sut$--representations.

\section{Hyperbolic structures on tori ($t < 2$)}

The topology and dynamics change dramatically as $t$ changes from
$t<-2$ to $t>2$. The level sets $t<-2$ correspond to Fricke spaces of
one-holed tori with geodesic boundary. The level set for $t=-2$
consists of four copies of the Teichm\"uller space of the punctured torus,
together with $\{(0,0,0)\}$. (The origin is the character of the
quaternion representation in $\sut$.) The origin is fixed under the 
action, while $\Gamma$ acts properly on its complement in $\kappa^{-1}(-2)$.
(However the action on the set of {\em complex points\/} of $\kappa^{-1}(-2)$
is extremely nontrivial and mysterious; see Bowditch~\cite{Bowditch}.)

The level sets for $-2<t<2$ correspond to Teichm\"uller spaces of singular
hyperbolic structures with one singularity, 
as well as a (compact) 
component consisting of characters of unitary representations. 
Except for the component of unitary 
representations, there are four components, freely permuted by $\Sigma$.
Except for the component of unitary representations, the $\Gamma$--action
is proper.

\subsection{Complete hyperbolic structures ($t\le -2$)}
The theory of deformations of geometric structures implies that 
$\Gamma$ acts properly on certain components of $\kt$. When $t < 2$, the
moduli space $\kt$ contains 4 contractible noncompact components,
freely permuted by $\Sigma$.  (When $-2 \le t < 2$, an additional
compact component corresponds to $\sut$--representations.) These
contractible components correspond to $\slr$--representations,
and each one identifies with the Teichm\"uller space of
$M$, with certain boundary conditions. Specifically, if $t<-2$, then
these components correspond to hyperbolic structures on $\o{int}(M)$
with geodesic boundary of length $2\cosh^{-1}(-t/2)$. For $t=-2$,
these components correspond to complete hyperbolic structures on
$\intt(M)$ and identify with the usual Teichm\"uller space
$\Teich$. For $-2<t<2$, these components correspond to singular
hyperbolic structures on a torus whose singularity is an isolated
point with cone angle 
\begin{equation*}
\theta = 2\cos^{-1}(-t/2). 
\end{equation*}

\subsection{Complete structures and proper actions}
We begin with the case $t<-2$. 
Then each component of $\ktR$ paramet\-rizes complete
hyperbolic structures on $\intt(M)$ with a closed geodesic parallel to
$\partial M$ having length $2\cosh^{-1}(-t/2)$. 
Thus the union 
\begin{equation*}
\bigcup_{t\le -2} \kt 
\end{equation*}
consists of four copies of the Fricke space of $M$. 
This space
contains equivalence classes of marked complete 
hyperbolic structures on $\o{int}(M)$, not necessarily of
finite area.
The properness of the action of the mapping class group
$\M$ on the Fricke space of $M$
implies properness of the $\Gamma$--action of $\kt$.
(For proof, see \S 2.2 of Abikoff~\cite{Abikoff}, 
Bers--Gardiner~\cite{BersG}
6.5.6 (page 156) of Buser~\cite{Buser}, 
6.3 of Imayoshi--Tanigawa~\cite{IT},
2.4.1 of Harvey~\cite{Harvey}, 
or \S 2.7 of Nag~\cite{Nag}.)

When $t=-2$, then $\kt$ has five connected components. 
It is the union of the single (singular) point $(0,0,0)$ and
four copies of the Teichm\"uller space $\Teich$ of $M$. Teichm\"uller
space $\Teich$ consists of equivalence classes of marked complete
finite-area hyperbolic structures on $M$. The  group $\Gamma$
acts properly on $\kt$ and fixes $(0,0,0)$.

\subsection{Singular hyperbolic structures on tori ($-2<t<2$)}

We next consider the case $2 > t > -2 $.
Here $\kt$ has five connected components. One component $C_t$ is compact
and consists of unitary characters, while the other four components
correspond to singular hyperbolic structures on a torus with a cone point.
$C_t$ is diffeomorphic to $S^2$, and the symplectic structure is a smooth
area form. $\Gamma$ acts ergodically on $C_t$. See \S 5 of Goldman~\cite{Erg}
for a detailed discussion.

Consider next the four noncompact components of $\kt$.  This case is
similar to the previous case, except that the ends of $M$ are replaced
by cone points on a torus.
Using the properness of the mapping class group action on
Teichm\"uller space, the action on these components of the relative
character variety remains proper. However, none of the corresponding
representations in $\hmg$ are discrete embeddings. Generically these
representations are isomorphisms onto dense subgroups of $\slr$.

We first show (Theorem~\ref{thm:ellipticCom}) 
that every representation in these components is a lift (to $\slr$) of the
holonomy representation of a singular hyperbolic structure.
In \S\ref{sec:properness}  we deduce that 
$\Gamma$ acts properly on the level sets $\kt$ where $-2<t<2$.

\subsection{Construction of hyperbolic structures on $T^2$ with one
cone point}

Let $\theta>0$ and let $\Cc_\theta$ denote the space of hyperbolic
structures on $T^2$ with a conical singularity of cone angle $\theta$.
By results of McOwen~\cite{McOwen} and Troyanov~\cite{Troyanov},
$\Cc_\theta$ identifies with the Teichm\"uller space 
$\Tei(\theta)$, the deformation space 
of conformal structures on $T^2$ singular at 
one point with cone angle $\theta$.

\begin{thm}\label{thm:ellipticCom}
Suppose $\xi,\eta\in\slr$ and 
$\gamma:=[\xi,\eta]$ is elliptic of rotation angle $\theta$.
Let $M$ denote a one-holed torus whose fundamental group
has the presentation \eqref{eq:pi1torus}.
Then there exists a singular hyperbolic structure on $\bar{M}$ 
with a singularity of cone angle $\theta$ having holonomy
representation $\rho$ defined by 
\begin{equation*}
\rho(X) = \xi,\  \rho(Y) = \eta,\  \rho(K) = \gamma. 
\end{equation*}
\end{thm}

The proof will be based on the following.

\begin{lemma}\label{lem:embeddedQ}
Let $p\in\Ht$ be a point fixed by $\gamma$ and consider the
points
\begin{align*}
p_4 &= p, \\
p_3 & = \eta^{-1} p, \\
p_2 & = \xi^{-1}  \eta^{-1} p, \\
p_1 & = \eta \xi^{-1}  \eta^{-1} p.
\end{align*}
The four points $p_1,p_2,p_3,p_4$ are the vertices of
an embedded quadrilateral $\qq$. In other words the four segments
\begin{align*}
l_1 &= \overline{p_1 p_2} \\
l_2 &= \overline{p_2 p_3} \\
l_3 &= \overline{p_3 p_4} \\
l_4 &= \overline{p_4 p_1}
\end{align*}
are disjoint and bound a quadrilateral.
\end{lemma}

\begin{proof}
[Proof of Theorem~\ref{thm:ellipticCom} assuming Lemma~\ref{lem:embeddedQ}]
\qquad

$\xi$ maps $p_1$ (respectively $p_2$) to $p_4$ (respectively $p_3$).
Therefore $\xi$ maps the directed edge $l_1$ to 
$l_3$ with the opposite orientation.
Similarly $\eta$ maps $p_2$ (respectively $p_3$) to $p_1$ (respectively $p_4$) 
so $\eta$ maps $l_2$ to $l_4$ with the opposite orientation.

From the embedding of $\qq$ and the identifications of its opposite
sides by $\xi$ and $\eta$, we construct a developing map for a (singular)
hyperbolic structure on $M$ as follows. 
Let $\qqq$ denote an {\em abstract quadrilateral,\/}
that is, a cell complex with a single 2--cell, four 1--cells
(the {\em sides\/}), and four 0--cells (the {\em vertices\/}).
Denote the oriented edges $\ll_i$ numbered in cyclic order 
and the vertices $p_i$ ($i=1,2,3,4$) where 
\begin{align*}
\partial \ll_1 = \{p_1,p_2\},& \quad \partial \ll_2 = \{p_2,p_3\}  \\
\partial \ll_3 = \{p_3,p_4\},& \quad \partial \ll_4 = \{p_4,p_1\}.  
\end{align*}
Let $\bar\ll_i$ denote $\ll_i$ with the opposite orientation. Write
\begin{equation*}
\qd = \qqq - \{p_1,p_2,p_3,p_4\}. 
\end{equation*}
Choose orientation-preserving homeomorphisms 
\begin{equation*}
\tilde\xi\co \ll_1 \longrightarrow \bar{\ll}_3, \qquad
\tilde\eta\co  \ll_2 \longrightarrow \bar{\ll}_4.
\end{equation*}
Let $\Pi=\langle X,Y\rangle$ denote the free group generated by $X,Y$. 
A model for the universal covering $\tilde M$ of $M$ 
is the quotient space of the Cartesian product
\begin{equation*}
\qd  \times \Pi
\end{equation*}
by the equivalence relation generated by the identifications
\begin{equation*}
(x,\eta) \longleftrightarrow (\tilde\xi(x),X\eta)
\end{equation*}
for $x\in \ll_1$, and
\begin{equation*}
(x,\eta) \longleftrightarrow (\tilde\eta(x),Y\eta)
\end{equation*}
for $x\in \ll_2$.  The action of $\Pi$ on $\qqq \times
\Pi$ defined by the trivial action on $\qqq$ and right-multiplication on 
$\Pi$ descends to a free proper action of $\Pi$ on $\tilde M$. 
This action corresponds to the action of deck transformations. 
Let $\rho\co \Pi\longrightarrow G$ be the homomorphism defined by:
\begin{align*}
\rho(X) & = \xi \\
\rho(Y) & = \eta 
\end{align*}
The embedding $\qd\longrightarrow\qq$ together with the
identifications of the sides of $\qq$ extends uniquely to a
$\rho$--equivariant local homeomorphism $\tilde M\longrightarrow
\Ht$. The holonomy around the puncture
\begin{equation*}
\rho([X,Y])=[\xi,\eta]=\gamma 
\end{equation*}
is elliptic of rotation angle $\theta$ fixing $p$. 
The resulting hyperbolic structure extends to a singular hyperbolic structure
on $\bar{M}$ with a singularity of cone angle $\theta$. 
Thus Theorem~\ref{thm:ellipticCom} 
reduces to Lemma~\ref{lem:embeddedQ}. \end{proof}

\begin{figure}[ht!]
\centerline{\epsfxsize=2.5in \epsfbox{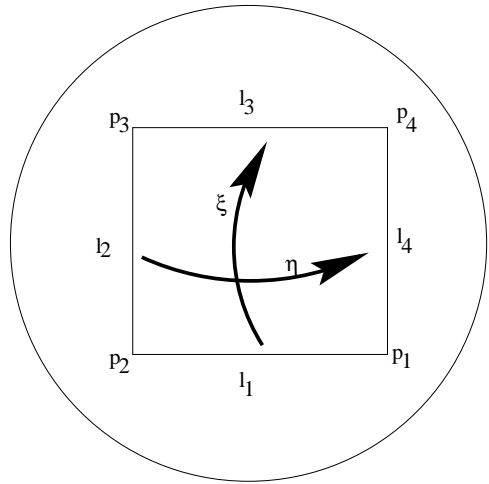}}
\caption{A convex embedded quadrilateral}
\label{fig:convex}
\end{figure}

\begin{figure}[ht!]
\centerline{\epsfxsize=2.5in \epsfbox{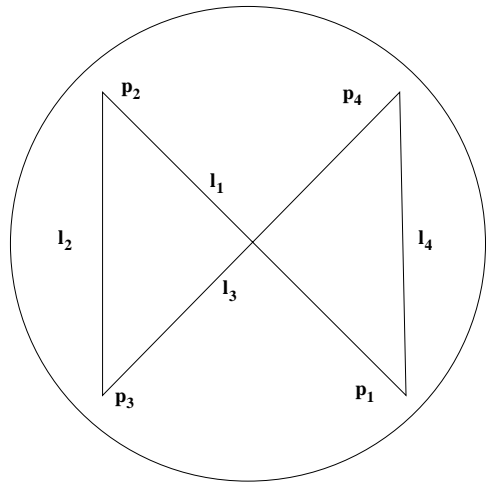}}
\caption{When $\overline{p_3p_4}$ meets $\overline{p_1p_2}$}
\label{fig:p34meetsp12}
\end{figure}

\begin{figure}[ht!]
\centerline{\epsfxsize=2.5in \epsfbox{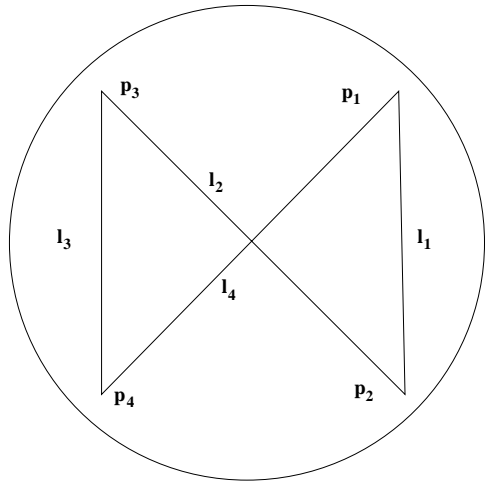}}
\caption{When $\overline{p_2p_3}$ meets $\overline{p_4p_1}$}
\label{fig:p23meetsp14}
\end{figure}

\begin{figure}[ht!]
\centerline{\epsfxsize=2.5in \epsfbox{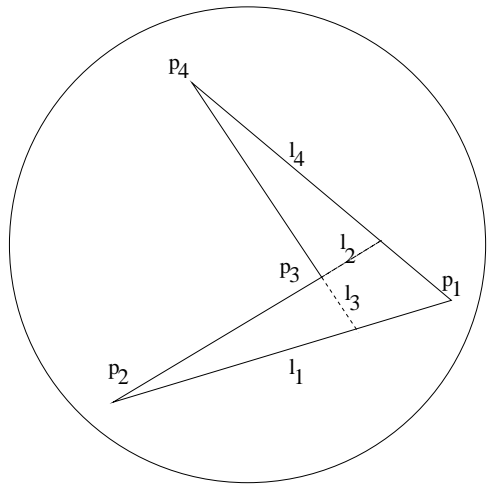}}
\caption{A nonconvex embedded quadrilateral}
\label{fig:nonconvex}
\end{figure}

The proof of Lemma~\ref{lem:embeddedQ} is based on the following two lemmas.

\begin{lemma}\label{lem:notcollinear}
The points $p_1,p_2,p_3,p_4$ are not collinear.
\end{lemma}

\begin{lemma}\label{lem:quadrilaterals}
If the sequence of points $p_1, p_2, p_3, p_4$ are not the 
vertices of an embedded quadrilateral, then either
$l_1$ and $l_3$ intersect, or
$l_2$ and $l_4$ intersect.
\end{lemma}
\begin{proof}[Proof of Lemma~\ref{lem:notcollinear}]
Suppose that $p_1,p_2,p_3,p_4$ are collinear.
If more than one line contains $p_1,p_2,p_3,p_4$, then all the points
coincide and $\rho$ fixes this point, contradicting $\gamma\neq\Id$.
Thus a unique line $l$ contains $p_1,p_2,p_3,p_4$. We claim that $l$
is $\rho$--invariant.

If $p_3 = p_4$, then $\eta$ fixes $p$ and since
$\gamma = \xi\eta\xi^{-1}\eta^{-1}$ fixes $p$, so must 
$\xi\eta\xi^{-1}$ fix $p$. Since $\gamma=[\xi,\eta]$
is nontrivial, both $\xi$ and $\eta$ are nontrivial. Thus $\xi$ maps the
unique fixed point of $\eta$ to the unique fixed point of $\xi\eta\xi^{-1}$,
and both $\xi$ and $\eta$ fix $p$, contradicting $\gamma\neq\Id$.
Thus $p_3\neq p_4$ and $l$ is the unique line containing $p_3$ and $p_4$.

$\xi(p_2)=p_3$ and $\xi(p_1)=p_4$ imply that $p_1\neq p_2$. Thus $l$
is the unique line containing $p_1$ and $p_2$, and is $\xi$--invariant.
Similarly $l$ is $\eta$--invariant and thus $\rho(\pi)$--invariant,
as claimed.

The commutator subgroup of the stabilizer of a line $l$ consists of
hyperbolic elements, contradicting ellipticity of $\gamma$.
\end{proof}

\begin{proof}[Proof of Lemma~\ref{lem:quadrilaterals}]
If $x,y\in\Ht$ are distinct points, let $\overleftrightarrow{xy}$
denote the line passing through $x$ and $y$. If
$z\notin\overleftrightarrow{xy}$, let $H_{x,y}(z)$ denote the
component (a half-space) of $\Ht-\overleftrightarrow{xy}$ containing
$z$, and $\hat{H}_{x,y}(z)$ denote the complementary component so that:
\begin{equation*}
\Ht =  H_{x,y}(z)\;  \amalg \; \overleftrightarrow{xy}\;  \amalg \; \hat{H}_{x,y}(z)
\end{equation*}
Suppose first that no three of the points $p_1,p_2,p_3,p_4$ are collinear.
There are four cases, depending on where $p_4$ lies in relation
to the decompositions determined by the lines $\overleftrightarrow{p_1p_2}$
and $\overleftrightarrow{p_2p_3}$:
\begin{enumerate}
\item $p_4\in H_{p_1p_2}(p_3)\cap H_{p_2p_3}(p_1)$;
\item $p_4\in \hat{H}_{p_1p_2}(p_3)\cap H_{p_2p_3}(p_1)$;
\item $p_4\in H_{p_1p_2}(p_3)\cap \hat{H}_{p_2p_3}(p_1)$;
\item $p_4\in \hat{H}_{p_1p_2}(p_3)\cap\hat{H}_{p_2p_3}(p_1)$.
\end{enumerate}
In the first case $p_1,p_2,p_3,p_4$ are the vertices of a convex quadrilateral.
(Compare Figure~\ref{fig:convex}.)
In the second case, $\overline{p_3p_4}$ meets $\overline{p_1p_2}$.
(Compare Figure~\ref{fig:p34meetsp12}.)
In the third case, $\overline{p_2p_3}$ meets $\overline{p_4p_1}$.
(Compare Figure~\ref{fig:p23meetsp14}.)
In the fourth case $p_1,p_2,p_3,p_4$ are the vertices of an embedded
(nonconvex) quadrilateral.
(Compare Figure~\ref{fig:nonconvex}.)

Suppose next that three of the vertices are collinear. By 
Lemma~\ref{lem:notcollinear}, not all four vertices are collinear.
By possibly conjugating $\rho$ by $\xi,\eta$ etc., we may assume that
$p_1,p_2,p_3$ are collinear. If $p_2$ lies between $p_1$ and $p_3$,
then $p_1,p_2,p_3,p_4$ are the vertices of an embedded quadrilateral (with
a straight angle at $p_2$). Otherwise
$\overline{p_2p_3}$ meets $\overline{p_4p_1}$ at $p_1$.
This completes the proof of Lemma~\ref{lem:quadrilaterals}.
\end{proof}
Returning to the proof of Lemma~\ref{lem:embeddedQ}, we 
show that $l_1$ and $l_3$ cannot intersect; 
an identical proof implies $l_2$ and $l_4$ cannot intersect.
\begin{claim*}
Suppose that $l_1$ and $l_3$ intersect. Then $l_1 \cap l_3$ is a point $q$ and
the triangles $\triangle(qp_1p_4)$ and
$\triangle(qp_3p_2)$ are congruent.
\end{claim*}
\begin{proof}[Proof of Claim]
If $l_1\cap l_3$ is not a point, then it is a segment, contradicting
Lemma~\ref{lem:notcollinear}. 
Since $\xi(l_1)=l_3$ and $\eta(l_2)=l_4$, 
the lengths of opposite sides are equal:
\begin{align*}
d(p_1,p_2) & = \ell(l_1) = \ell(l_3) = d(p_4,p_3), \\
d(p_2,p_3) & = \ell(l_2) = \ell(l_4) = d(p_1,p_4)
\end{align*}
Since the lengths of the corresponding sides are equal,
the triangle $\triangle(p_2p_1p_4)$ is congruent to
$\triangle(p_4p_3p_2)$. In particular
\begin{equation*}
\angle(p_2p_1p_4) = \angle(p_4p_3p_2). 
\end{equation*}
Similarly, $\angle(p_3p_4p_1) = \angle(p_1p_2p_3)$.
Now 
\begin{equation*}
\angle(q p_1p_4) = \angle(p_2p_1p_4) = 
\angle(p_4p_3p_2) = \angle(qp_3p_2) 
\end{equation*}
and similarly, $\angle(q p_4p_1) = \angle(qp_2p_3)$. 
Since $l_4 = \overline{p_1p_4}$ is congruent to 
$l_2 =\overline{p_3p_2}$, triangles
$\triangle(qp_1p_4)$ and
$\triangle(qp_3p_2)$ are congruent as claimed.
\end{proof}

Now let $l$ be the line through $q$ bisecting the angle
$\angle(p_1qp_3)$ such that reflection $R$ in $l$
interchanges $\triangle(qp_1p_4)$ and $\triangle(qp_3p_2)$.
Since 
\begin{equation*}
R\co\begin{cases}
p_1 &\longmapsto  p_3 \\
p_2 &\longmapsto  p_4 \\
p_3 &\longmapsto  p_1 \\
p_4 &\longmapsto  p_2, 
\end{cases}
\end{equation*}
\begin{equation*}
\xi\co  
\begin{cases}
p_1 &\longmapsto  p_4 \\
p_2 &\longmapsto  p_3
\end{cases}
\end{equation*}
\begin{equation*}
\eta\co 
\begin{cases}
p_2 &\longmapsto  p_1 \\
p_3 &\longmapsto  p_4,
\end{cases}
\end{equation*}
$\xi\circ R$ interchanges $p_3$ and $p_4$, and
$R\circ\eta$ interchanges $p_2$ and $p_3$. 
Since an orientation-reversing isometry of $\Ht$
which interchanges two points must be reflection
in a line, $\xi\circ R$ and 
$R\circ\eta$ have order two. Thus
$R$ conjugates $\xi$ to $\xi^{-1}$ and 
$\eta$ to $\eta^{-1}$.

One of two possibilities must occur:
\begin{itemize}
\item At least one of $\xi$ and $\eta$ is elliptic or parabolic;
\item Both $\xi$ and $\eta$ are hyperbolic, and their invariant axes
are each orthogonal to $l$.
\end{itemize}
Neither possibility occurs, due to the following:

\begin{lemma}\label{lem:crossing axes}
Let $\xi,\eta\in\slr$. The following conditions are equivalent:
\begin{itemize} 
\item 
$\tr [\xi,\eta] < 2$;
\item
$\xi,\eta$ are hyperbolic elements and their invariant axes cross.
\end{itemize} 
\end{lemma}
\begin{proof}
Assuming $\tr [\xi,\eta] < 2$, we first show that $\xi$ and $\eta$
must be hyperbolic. We first show that $\xi$ is not elliptic.
If $\xi$ is elliptic (or $\pm\Id)$, 
we may assume that $\xi\in\sot$, 
that is, we represent $\xi,\eta$ by matrices
\begin{equation*}
\xi = \bmatrix \cos(\theta) & -\sin(\theta) \\ 
\sin(\theta) & \cos(\theta) \endbmatrix,\quad
\eta = \bmatrix a  & b \\ c  & d \endbmatrix,
\end{equation*}
where $ad-bc = 1 $, whence
\begin{equation*}
\tr [\xi,\eta] = 2 + \sin^2(\theta) (a^2 + b^2 + c^2 + d^2 - 2) \ge 2. 
\end{equation*}
Similarly if $\xi$ is parabolic, \eqref{eq:trcomunip} implies that
$\tr[\xi,\eta]\ge 2$.

Thus $\xi$ is hyperbolic.  Since $[\eta,\xi] =
[\xi,\eta]^{-1}$, an identical argument shows that $\eta$ is
hyperbolic.  Denote their invariant axes by $l_\xi,l_\eta$
respectively.

It remains to show that $\tr[\xi,\eta]<2$ if and only if
$l_\xi\cap l_\eta\neq\emptyset$.

By conjugation, we may assume that the fixed points of $\xi$ are $\pm 1$
and that the fixed points of $\eta$ are $r,\infty$. Thus
$l_\xi\cap l_\eta\neq\emptyset$  if and only if $-1< r <1$.
Represent $\xi,\eta$ by matrices
\begin{equation*}
\xi = \bmatrix \cosh(\theta) & \sinh(\theta) \\ 
\sinh(\theta) & \cosh(\theta) \endbmatrix, \quad
\eta = \bmatrix e^\phi  & -2r \sinh(\phi) \\ 0  & e^{-\phi} \endbmatrix,
\end{equation*}
where $\theta,\phi\neq 0$ and
\begin{equation*}
\tr [\xi,\eta] = 2 + 4 (r^2-1) \sinh^2(\theta) \sinh^2(\phi).
\end{equation*}
Then $\tr [\xi,\eta] < 2$ if and only if $-1<r<1$ and 
$\tr [\xi,\eta] = 2$ if and only if $r=\pm 1$, as desired.
\end{proof}

This completes the proof of Lemma~\ref{lem:embeddedQ} 
(and also Theorem~\ref{thm:ellipticCom}).\qed

\subsection{Properness}\label{sec:properness}

That $\Gamma$ acts properly on $\ktR$ now follows easily. 
The construction in Theorem~\ref{thm:ellipticCom} gives a map
\begin{equation}\label{eq:uniformization}
\kt \cap[2,\infty)^3  \longrightarrow \Tei(\theta)
\end{equation}
which is evidently $\Ou$--equivariant.
Every hyperbolic structure with conical singularities has an underlying
{\em singular conformal structure,\/} where the singularities are
again conical singularities, that is, they are defined by local coordinate
charts to a model space, which in this case is a cone. 
However, there is an important difference. 
Conical singularities in conformal structures are {\em removable,\/}
while conical singularities in Riemannian metrics are not.

Here is why conformal conical singularities are removable:
Let $D^2\subset\C$ be the unit disk and let 
\begin{equation*}
S_{\theta} := \{ z\in D^2 \mid 0 \le \arg(z)\le \theta\}
\end{equation*}
be a sector of angle $\theta$. Then
\begin{align*}
\Pi_\theta:S_{\theta} - \{0\} & \longrightarrow D^2- \{0\}  \\
                             z &\longmapsto z^{2\pi/\theta}
\end{align*}
is conformal. 
The model coordinate patch for a cone point of angle $\theta$
is the  {\em cone\/} $C_\theta$ 
of angle $\theta$, defined as the identification space 
$C_{\theta}$ of $S_{\theta}$ by the equivalence relation defined by 
\begin{equation*}
z \longleftrightarrow e^{\pm i\theta} z 
\end{equation*}
for $z\in\partial S_\theta$.
That is, a cone point $p$ has a coordinate patch neighborhood $U$ and
a coordinate chart $\psi\co U\longrightarrow C_{\theta}$ in the atlas
defining the singular geometric structure.
The power map $\Pi_\theta$ defines a conformal isomorphism between the 
$C_{\theta}$ punctured at the cone point and a punctured disk (a ``cone'' of
angle $2\pi$).
Replacing the  coordinate chart $\psi\co U\longrightarrow \C$ 
at a cone point $p$ of angle $\theta$ by the composition $\Pi_\theta\circ\psi$
gives a coordinate atlas for a conformal structure which is nonsingular at
$p$ and isomorphic to the original structure on the complement of $p$.

The resulting map $\Tei(\theta) \longrightarrow \Teich$
is evidently $\Ou$--equivariant. 
Since $\Ou$ acts properly on $\Teich$, 
$\Ou$ acts properly on $\Tei(\theta)$, and hence on $\ktR$ as well.

With more work one can show that \eqref{eq:uniformization} is an
isomorphism. For any hyperbolic structure on $T^2$ with a cone point
$p$ of angle $0<\theta<2\pi$, the
Arzel\'a--Ascoli theorem (as in Buser~\cite{Buser} \S 1.5) applies
to represent $\xi,\eta$ by geodesic loops based at
$p$. Furthermore these geodesics intersect only at $p$. By developing
this singular geometric structure one obtains a polygon $\qq$ as above.

\section{Reducible characters ($t=2$)}
The level set $\ktwo$ 
consists of characters of reducible representations.
Over $\C$ such a representation is upper-triangular
\eqref{eq:uppertri} with character defined by \eqref{eq:extension}.
By \eqref{eq:factorization}, every $(x,y,z)\in\ktwo$ lies
in the image of the map
\begin{align}\label{eq:reduciblechar}
\Phi\co  \C^* \times \C^* & \longrightarrow \ktwo \\
(\xi,\eta) & \longmapsto 
\bmatrix \xi + \xi^{-1} \\ \eta + \eta^{-1} \\ 
\xi\eta + \xi^{-1}\eta^{-1} \endbmatrix. \notag
\end{align}
The set of $\R$--points $\ktwoR$ is a singular algebraic
hypersurface in $\R^3$, with singular set
\begin{equation*}
S_0 =\left\{\bmatrix 2 \\ 2 \\2\endbmatrix, 
       \bmatrix 2 \\ -2 \\-2\endbmatrix, 
       \bmatrix -2 \\ 2 \\ -2 \endbmatrix,
       \bmatrix -2 \\ -2\\ 2\endbmatrix \right\}.
\end{equation*}
Characters in $S_0$
correspond to unipotent representations twisted by central characters
(as in the sense of \eqref{eq:twisted}).
A {\em central character\/} is a homomorphism taking values in 
the center $\{\pm\Id\}$ of $G$ and a {\em unipotent representation\/} is
a representation in $\pm U$, where $U$ is a unipotent subgroup of $G$.
A reductive representation with character in $S_0$ is itself a central
character. The most general representation with character in $S_0$
is one taking values in $\pm U$, where $U$ is a unipotent subgroup of $G$.
The character $(2,2,2)$ is the character of any unipotent representation, for
example the trivial representation.  The other three points are images
of $(2,2,2)$ by the three nontrivial elements of $\Sigma$. 

The smooth stratum of $\ktwo$ is the complement 
\begin{equation*}
\ktwoR - S_0. 
\end{equation*}
Denote its five components
by $C_K$ and $C_i$, where $i=0,1,2,3$. Here $C_0$ denotes the
component $\ktwo\cap [2,\infty)^3$ and $C_i = (\sigma_i)_* C_0$ where
$(\sigma_i)_*\in\Sigma$ is the sign-change automorphism fixing
the $i$-th coordinate (\S\ref{sec:signchanges}).

The component $C_K$ corresponds to reducible $\sut$--representations which
are {\em non-central,\/} that is, their image does not lie in the
center $\{\pm\Id\}$ of $\sut$. The closure of $C_K$ is the union of $C_K$
with $S_0$. The map
\begin{align*}
\uo \times \uo & \longrightarrow \bar{C}_K  \\
(\xi,\eta) & \longmapsto (\xi + \xi^{-1}, \eta + \eta^{-1},
\xi\eta + \xi^{-1}\eta^{-1})
\end{align*}
is a double branched covering space, with deck transformation
\begin{equation*}
(\xi,\eta) \longmapsto (\xi^{-1}, \eta^{-1})
\end{equation*}
and four branch points 
\begin{equation*}
(\xi,\eta) = (\pm 1,\pm 1)
\end{equation*}
which map to $S_0$.

\begin{figure}[ht!]
\centerline{\epsfxsize=3in \epsfbox{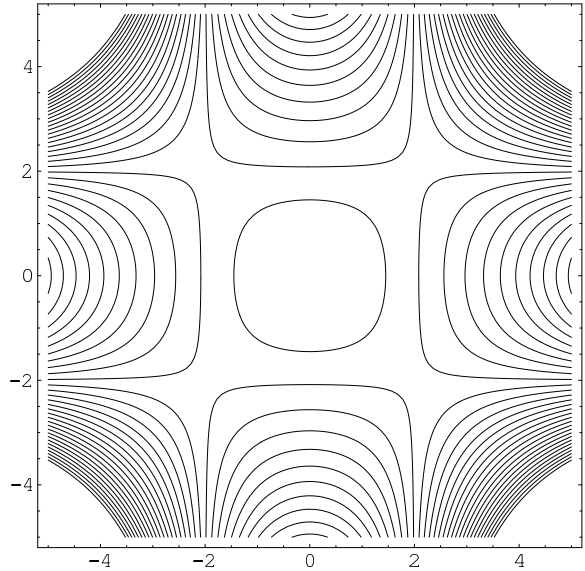}}
\caption{Level sets of $(x^2-4)(y^2-4)$}
\label{fig:contours}
\end{figure}

Similarly, each of the other components identifies with the quotient
$(\R_+)^2/\{\pm\Id\}$ with the action of $\GLtz$. For example,
$\bar{C}_0$ is the image of the double branched covering
\begin{align*}
\R_+ \times \R_+ & \longrightarrow \bar{C}_0 \\
(\xi,\eta) & \longmapsto (\xi + \xi^{-1}, \eta + \eta^{-1},
\xi\eta + \xi^{-1}\eta^{-1}),
\end{align*}
with deck transformation
\begin{equation*}
(\xi,\eta) \longmapsto (\xi^{-1}, \eta^{-1}).
\end{equation*}
For $i=1,2,3$, the component $C_i$ is the image of $C_0$ under the
sign-change $(\sigma_i)_*$.
Composing with the 
$\GLtz$--equivariant diffeomorphisms
\begin{align*}
\R^2 & \longrightarrow \R_+ \times \R_+ \\
(\tilde\xi,\tilde\eta) & \longmapsto 
(\exp(\tilde\xi),\exp(\tilde\eta))
\end{align*}
and
\begin{align*}
(\R/\Z)^2 & \longrightarrow \uo \times \uo \\
(\tilde\xi,\tilde\eta) & \longmapsto 
(\exp(2\pi i\tilde\xi),\exp(2\pi i\tilde\eta))
\end{align*}
respectively,
yields  $\GLtz$--equivariant double branched coverings
\begin{equation*}
\R^2\longrightarrow C_i, \text{~and~} (\R/\Z)^2\longrightarrow C_K  
\end{equation*}
respectively. Furthermore these
mappings pull back the invariant area form on $\ktwoR$ 
to Lebesgue measure on $\R^2$.

Since $\SL{2,\Z}$ is a lattice in $\slr$ and $\slr$ acts 
transitively on $\R^2$ with noncompact isotropy group,
Moore's ergodicity theorem (Moore~\cite{Moore}; 
see also Feres~\cite{Feres}, Zimmer~\cite{Zimmer} or Margulis~\cite{Margulis}) 
implies that $\SL{2,\Z}$ acts ergodically on $\R^2$.
Thus $\GLtz$ acts ergodically on $\R^2/\{\pm\Id\}$ and
$(\R/\Z)^2/\{\pm\Id\}$, and hence on each of the components 
\begin{equation*}
C_0,C_1,C_2,C_3,C_K \subset \ktwoR. 
\end{equation*}
Since $\Sigma$ permutes $C_0,C_1,C_2,C_3$, 
the $\Gamma$--action on their union
\begin{align*}
C_0 \cup C_1 \cup C_2 \cup C_3 & = \ktwoR - C_K \cup S_0 \\ & 
\cong  (\R^*\times\R^*)/\{\pm\Id\}
\end{align*}
is ergodic. 
(In the case of $C_K$, any hyperbolic element in $\GLtz$ acts 
ergodically on $\uo\times\uo$ and hence on $\bar{C}_0$, a much stronger
result.)

\section{Three-holed spheres and ergodicity ($t>2$)}
Next we consider the level sets where $t>2$. There is an important
difference between the cases when $t> 18$ and $2<t\le 18$.  When $t
\le 18$, the $\Gamma$--action is ergodic, but when $t>18$, wandering
domains appear, arising from the Fricke spaces of a three-holed sphere
$P$ (``pair-of-pants'').  The three-holed sphere is the only other
orientable surface homotopy-equivalent to a one-holed torus, and
homotopy equivalences to hyperbolic manifolds homeomorphic to $P$
define points in these level sets when $t> 18$, which we call {\em
discrete $P$--characters\/} However, the $\Gamma$--action on the
complement of the discrete $P$--characters is ergodic.

\subsection{The Fricke space of a three-holed sphere.}
When $t>18$, the octant 
\begin{equation*}
\Omega_0 = (-\infty,-2)\times(-\infty,-2)\times(-\infty,-2)  
\end{equation*}
intersects $\kt$ in a wandering domain and the images of $\Omega_0\cap\kt$ are
freely permuted by $\Gamma$. Let
\begin{equation*}
\Omega = \Gamma\cdot\Omega_0.
\end{equation*}
Characters in $\Omega$ correspond to discrete embeddings 
$\rho\co \pi\longrightarrow\slr$ where the complete hyperbolic surface
$\Ht/\rho(\pi)$ is diffeomorphic to a three-holed sphere.
We call such a discrete embedding a {\em discrete $P$--embedding,\/}
and its character a {\em discrete $P$--character.\/}

The fundamental group $\pi_1(P)$ is free of rank two. A pair of
boundary components $\partial_1,\partial_2$, an orientation on $P$,
and a choice of arcs $\alpha_1,\alpha_2$ from the basepoint to  
$\partial_1,\partial_2$ determines a pair of free generators of $\pi_1(P)$:
\begin{align*}
X &:= (\alpha_1)^{-1} \star \partial_1 \star \alpha_1 \\
Y &:= (\alpha_2)^{-1} \star \partial_2 \star \alpha_2.
\end{align*}
A third generator $Z := (XY)^{-1}$ corresponds to the third boundary
component, obtaining a presentation of $\pi_1(P)$ as
\begin{equation*}
\langle X,Y,Z \mid XYZ = 1 \rangle. 
\end{equation*}
Elements of $\Omega_0$ are discrete $P$--characters such that the generators
$X,Y$ and $Z := (XY)^{-1}$ of $\pi$ correspond to the boundary components of
the quotient hyperbolic surface $\Ht/\rho(\pi)$. 

\begin{figure}[ht!]
\centerline{\epsfxsize=2.5in \epsfbox{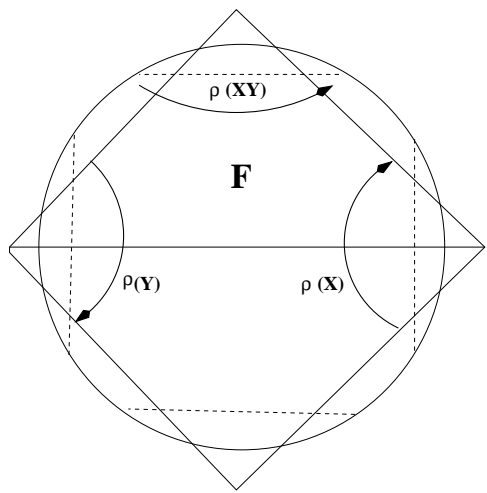}}
\caption{Fundamental hexagons for a hyperbolic three-holed sphere}
\label{fig:hexagon}
\end{figure}

\begin{lemma}\label{lem:pantschar} 
A representation $\rho\in\Hm{\slt}$ has character
$[\rho] = (x,y,z)\in\Gamma\cdot\Omega_0$ if and only if
$\rho$ is a discrete $P$--embedding such that $X,Y,Z$ correspond to
the components of $\partial\Ht/\rho(\pi)$.
\end{lemma}
\begin{proof}
The condition that 
$[\rho]\in\Omega_0$ is equivalent to $x,y,z < -2$,  which implies that
the generators $\rho(X)$, $\rho(Y)$ and $\rho(XY)$ are hyperbolic and
their invariant axes $l_{\rho(X)}, l_{\rho(Y)}, l_{\rho(XY)}$  
are pairwise ultraparallel. Denoting the common perpendicular to
two ultraparallel lines $l,l'$ by $\perp(l,l')$, the six lines
\begin{equation*}
l_{\rho(X)},\thickspace 
\perp\negthickspace\big(l_{\rho(X)},l_{\rho(Y)}\big), \thickspace
l_{\rho(Y)},\thickspace 
\perp\negthickspace\big(l_{\rho(Y)},l_{\rho(XY)}\big),  \thickspace
l_{\rho(XY)},\thickspace  
\perp\negthickspace\big(l_{\rho(XY)},l_{\rho(X)}\big) 
\end{equation*}
bound a right-angled hexagon $H$. 
The union of $H$ with its reflected image in 
$\perp\big(l_{\rho(X)},l_{\rho(Y)}\big)$ is a fundamental domain for
$\rho(\pi)$ acting on $\Ht$.  
(Goldman~\cite{Gth}, Gilman--Maskit~\cite{GM}, I-7, page 15). 
(Figure~\ref{fig:hexagon} depicts the identifications corresponding to
the generators $\rho(X),\rho(Y)$.)  
The quotient is necessarily homeomorphic to a three-holed sphere $P$ 
and the holonomy around components of $\partial P$ are the
three generators $\rho(X),\rho(Y),\rho(XY)$.
\end{proof}

Recall (\S\ref{sec:permutations}) that $\Gamma$ decomposes as the
semidirect product
\begin{equation*}
\Gamma \cong \Gamma_{(2)} \ltimes \Ss_3.
\end{equation*}
The mapping class group of a three-holed sphere 
$P$ is isomorphic to $\Z/2\times\Ss_3$.
The $\Ss_3$--factor corresponds to the group of permutations of 
$\pi_0(\partial P)$. The $\Z/2$--factor is generated by the elliptic
involution 
which acts by an orientation-reversing
homeomorphism of $P$, whose fixed-point set is the union of three
disjoint arcs joining the boundary components.
(See \S\ref{sec:other} for the corresponding automorphism of $\pi$.) 
Accordingly, $\Ss_3$ preserves the Teichm\"uller space of $P$. 
The following result shows that these are the only automorphisms
of the character variety which preserve the discrete $P$--characters.

\begin{proposition} $\Omega$ equals the disjoint union 
$\coprod_{\gamma\in\Gamma_{(2)}}\gamma\Omega_0$.
\end{proposition}
\begin{proof}
We show that if for some $\gamma\in\Gamma_{(2)}$, the intersection
$\Omega_0\cap \gamma\Omega_0$ is nonempty, then $\gamma=1$.  Suppose
that $[\rho]\in\Omega_0\cap \gamma\Omega_0$. By Lemma~\ref{lem:pantschar},
both $\rho$ and $\rho\circ\gamma$ are discrete $P$--embeddings such that
$X,Y,(XY)^{-1}$ and $\gamma(X),\gamma(Y),\gamma((XY)^{-1})$ correspond
to $\partial\Ht/\rho(\pi)$.

The automorphism $\gamma$ of the character space $\C^3$ corresponds to
an automorphism $\tg$ of $\pi$ such that 
\begin{equation*}
[\rho\circ\tg] = \sigma_*\circ\gamma([\rho])   
\end{equation*}
for a sign-change automorphism $\sigma_*\in\Sigma$.
Thus $\rho\circ\tg$ is also a discrete $P$--embedding, with quotient
bounded by curves corresponding to $\tg(X), \tg(Y)$ and $\tg(XY)$. 
Then $\tg(X)$ (respectively $\tg(Y)$,$\tg(XY)$) is conjugate to
$X^{\pm 1}$ (respectively $Y^{\pm 1}$, $(XY)^{\pm 1}$). 
Such an automorphism is induced by a diffeomorphism of the three-holed sphere 
$P$. Since the elliptic involution generates the 
mapping class group of $P$, 
$\tg$ must be an inner automorphism possibly composed with the
elliptic involution of $\pi$, and hence must act trivially on characters.
Thus $\gamma=1$ as desired.
\end{proof}

For $t\le 18$, the domain $\Omega$ does not meet $\kt$. (The closure
$\bar{\Omega}$ intersects $\kappa^{-1}(18)$ in the $\Gamma$--orbit
of the character $(-2,-2,-2)$ of the holonomy representation of a complete
finite-area hyperbolic structure on $P$.)
For $t>18$, observe that $\kt -\Omega$ contains the open subset
\begin{equation*}
\bigg((-2,2)\times\R\times\R\bigg)\cap\kt.  
\end{equation*}
Thus $\kt -\Omega$ has positive measure in $\kt$.

\begin{proposition}\label{prop:erg}
For any $t>2$, the action of $\Gamma$ on $\kt-\Omega$ is ergodic.
\end{proposition}

The proof uses an iterative procedure
(Theorem~\ref{thm:IterativeProcedure}) due to Kern-Isberner and
Rosenberger~\cite{K-IR}, although their proof contains a
gap near the end.  See also Gilman and Maskit~\cite{GM}.  
Theorem~\ref{thm:IterativeProcedure} is proved at the end of the paper.

\subsection{The equivalence relation defined by $\Gamma$}

Write $u\sim v$ if there exists $\gamma \in\Gamma $ such that $\gamma u = v$.
Since $\kappa$ is $\Gamma$--invariant, $u\sim v$ implies that
$\kappa(u)=\kappa(v)$.

\begin{thm}\label{thm:IterativeProcedure}
 Suppose that $u\in\R^3$ satisfies $\kappa(u)>2$. Then
there exists $(\tx,\ty,\tz) = \tu\sim u$ such that either
\begin{itemize}
\item $\tu\in \left(-\infty,-2\right]^3$, in which case $\tu$ is the
character of a Fuchsian representation whose quotient is a hyperbolic
structure on a three-holed sphere $P$, with boundary mapping to either cusps
or closed geodesics;
\item $\tx\in \left[-2,2\right]$, in which case $\tu$ is the character of
a representation mapping $X$ to a non-hyperbolic element.
\end{itemize}
\end{thm}

Recall that a measurable equivalence relation is {\em ergodic\/}
if the only invariant measurable sets which are unions of equivalence
classes are either null or conull. Equivalently an equivalence relation
is ergodic if and only if every function constant on equivalence classes
is constant almost everywhere. A group action defines an equivalence 
relation. However, equivalence relations are more flexible since every
subset of a space with an equivalence relation inherits an equivalence relation
(whether it is invariant or not). Suppose that $S$ is a measurable subset
of a measure space $X$.
If every point in $X$ is equivalent to
a point of a measurable subset of $S$, 
then ergodicity of the equivalence relation on $X$ 
is equivalent to ergodicity of $S$ with respect to the measure
class induced from $X$.

\begin{proof}[Proof of Proposition~\ref{prop:erg} assuming
Theorem~\ref{thm:IterativeProcedure}]
For $i=1,2$, let $\E^{(i)}_t$ denote the subset of $\kt-\Omega$ where at least 
$i$ of the coordinates lie in $[-2,2]$.
Suppose that $t>2$ and $u\in\big(\kt-\Omega\big)$.
Then Theorem~\ref{thm:IterativeProcedure} implies that 
$\Gamma u \cap \E_t^{(1)} \neq \emptyset$.
Thus ergodicity of the $\Gamma$--action on $\kt-\Omega$ is equivalent
to the ergodicity of the induced equivalence relation on
$\E_t^{(1)}$.

We now reduce to the equivalence relation on $\E_t^{(2)}$.
By applying a permutation we may assume that $-2 < x < 2$.

The level set 
\begin{equation*}
E(x_0) := \kt \cap x^{-1}(x_0)
\end{equation*}
is defined by
\begin{equation*}
\frac{2-x_0}4 (y + z)^2 +
\frac{2+x_0}4 (y - z)^2 = t - 2 + x_0^2
\end{equation*}
and is an ellipse since $-2 < x_0 < 2$. 
Furthermore the symplectic measure on $\kt$ disintegrates under
the map $x\co \kt \longrightarrow \R$ to $\tau_X$--invariant 
Lebesgue measure on $E(x_0)$.
In particular the Dehn twist $\tau_X$ (see \S\ref{sec:other}) 
acts by the linear map
\begin{equation*}
\bmatrix y \\ z \endbmatrix \longrightarrow
\bmatrix x & -1 \\ 1 & 0 \endbmatrix 
\bmatrix y \\ z \endbmatrix,
\end{equation*}
which is linearly conjugate to a rotation of the circle through 
angle $\cos^{-1}(x/2)$. Thus for almost every $x_0\in(-2,2)$ the Dehn
twist $\tau_X$ acts on $E(x_0)$ by a rotation of infinite order, 
and the action is ergodic. Furthermore by applying powers of $\tau_X$,
we may assume that $-2 < y < 2$ as well.
Thus ergodicity of the $\Gamma$--action on $\kt-\Omega$ is equivalent
to the ergodicity of the induced equivalence relation on
$\E_t^{(2)}$.

 Since $\Gamma$ acts by polynomial
transformations over $\Z$, those points of $\kt$ which are equivalent
to a point with $\cos^{-1}(x/2)$ rational comprise a set of measure zero.
We henceforth restrict to the complement of this set.

The quadratic reflection (see \S\ref{sec:other})
\begin{equation*}
Q_z\co  \bmatrix x \\ y \\ z \endbmatrix \longmapsto
\bmatrix x \\ y \\ x y - z \endbmatrix 
\end{equation*}
is the deck transformation for the double covering of $\kt$ given
by projection $\Pi_{(x,y)}$ to the $(x,y)$--plane. The image
\begin{equation*}
\Pi_{(x,y)}\co  \kt\subset\R^3 \longrightarrow \R^2 
\end{equation*}
is the region 
\begin{equation*}
\rr_t := \big\{(x,y)\in\R^2 \mid (x^2 - 4)(y^2 - 4) + t - 2 \ge 0\big\}
\end{equation*}
and $\Pi_{(x,y)}\co \kt\longrightarrow\rr_t$ is the quotient map for the action
of $Q_z$.
Thus ergodicity on $\E_t^{(2)}$ reduces to ergodicity of the induced
equivalence relation on
\begin{equation*}
\Pi_{(x,y)}(\E_t^{(2)}) = [-2,2]\times [-2,2].
\end{equation*}
Ergodicity now follows as in \S 5.2 of \cite{Erg}.  Suppose that
$f\co \kt-\Omega\longrightarrow\R$ is a $\Gamma$--invariant measurable
function. The ergodic decomposition for the equivalence relation
induced by the cyclic group $\langle \tau_X\rangle$ is the coordinate
function 
\begin{equation*}
x\co [-2,2]\times [-2,2]\longrightarrow [-2,2],  
\end{equation*}
and by ergodicity of $\langle\tau_X\rangle$ on the level sets of $x$,
there is a measurable function $g\co [-2,2]\longrightarrow \R$ such that
$f$ factors as $f= g\circ x$ almost everywhere. Applying the cyclic
group $\langle\tau_Y\rangle$ to $[-2,2]$, the function $g$ is constant
almost everywhere. Hence $f$ is constant almost everywhere, proving
ergodicity.
\end{proof}

\subsection{The trace-reduction algorithm}

The proof of Theorem~\ref{thm:IterativeProcedure}
is based on the following:
\begin{lemma}\label{lem:MainLemma}
Let $2<x\le y\le z$. Suppose that $\kappa(x,y,z) > 2$. Let $z' = xy -z$.
Then $z-z' > 2\sqrt{\kappa(x,y,z) -2}$.
\end{lemma}
\noindent The following expression for $\kappa(x,y,z)$ will be useful:
\begin{align}\label{eq:kappaxy}
\kappa(x,y,z) - 2  
	& = \frac14 \left\{ (2z - xy)^2 - (x^2-4)(y^2-4)\right\}  \\
	& = \frac14 \left\{ (z - z')^2 - (x^2-4)(y^2-4)\right\}. \notag
\end{align}
For fixed $x,y>2$, write $\kappa_{x,y}\co \R\to\R$ for the quadratic function
\begin{equation*}
\kappa_{x,y}\co  z \longmapsto \kappa(x,y,z) -2 = z^2 -  x y z + (x^2 + y^2 - 4). 
\end{equation*}
Then 
\begin{equation*}
\kappa_{x,y}^{-1}\bigg((-\infty,0])\bigg) = 
[\zeta_-(x,y), \zeta_+(x,y)]
\end{equation*}
where 
\begin{equation*}
\zeta_{\pm}(x,y) = \frac12 \left( xy \pm \sqrt{(x^2-4)(y^2-4)}\right).
\end{equation*}
Furthermore reflection
\begin{equation*}
z\longmapsto z' = xy - z
\end{equation*}
interchanges the two intervals
\begin{align*}
J_-  & = \big( -\infty,\,\zeta_-(x,y) \big) \\
J_+  & = \big( \zeta_+(x,y),\,\infty \big)
\end{align*}
comprising $\kappa_{x,y}^{-1}\big((0,\infty)\big)$.
\begin{lemma}\label{lem:yinterval}
Suppose that $2<x\le y$. Then $\zeta_-(x,y)< y < \zeta_+(x,y)$
\end{lemma}
\proof
The conclusion is equivalent to $\kappa(x,y,y)<2$, which is what we prove.
First observe that $y> 2$ implies $y - 1/2  > 3/2$ so that
\begin{align*}
y - y^2 & = -\left(y-\frac12\right)^2 + \frac14 \\
& < -\left(\frac32\right)^2 + \frac14  = -2
\end{align*}
and $x<y$ implies that $x + 2 - y^2  < y + 2 - y^2 < 0$. 
Therefore
\begin{align*}
\kappa(x,y,y) - 2 & = x^2 + 2y^2 - xy^2 -4 \\
& = (x^2 -4) + (2y^2 - xy^2) \\
& = (x-2) (x + 2 - y^2) < 0.\tag*{\qed}
\end{align*}

\begin{proof}[Conclusion of Proof of Lemma~\ref{lem:MainLemma}]
$\phantom{99}$

By Lemma~\ref{lem:yinterval}, 
the quadratic function $\kappa_{x,y}$ is negative at $y$.
By hypothesis $\kappa_{x,y}$ is positive at $z\ge y$. 
Therefore 
\begin{equation*}
z > \zeta_+(x,y) > y.  
\end{equation*}
Thus $z\in J_+$. 
Reflection $z\longmapsto z'$ interchanges the intervals $J_+$ and $J_-$, 
so $z'\in J_-$, that is $z'<\zeta_-(x,y) < y$ (Lemma~\ref{lem:yinterval}).
By \eqref{eq:kappaxy},
\begin{align*}
(z - z')^2 & = 4  \bigg(\kappa(x,y,z) -2\bigg) + (x^2 - 4) (y^2 -4 ) \\
& > 4  (\kappa(x,y,z) -2)
\end{align*}
whence (because $z \ge y > z'$)
\begin{equation*}
z - z'  > 2 \sqrt{\kappa(x,y,z) -2},
\end{equation*}
completing the proof of Lemma~\ref{lem:MainLemma}. \end{proof}

\begin{proof}[Conclusion of Proof of Theorem~\ref{thm:IterativeProcedure}]
Fix 
\begin{equation*}
u=(x,y,z)\in\R^3 
\end{equation*}
with $\kappa(u)>2$. We seek $\tu\sim u$ such that 
one of the following possibilities occurs:

\begin{itemize}
\item one of the coordinates $\tx,\ty,\tz$ lies in the interval $[-2,2]$; 
\item $\tx,\ty,\tz < -2$.
\end{itemize}
It therefore suffices to find $\tilde u\sim u $ which lies in
$(-\infty,2]^3$. Suppose that $u$ does not satisfy this. 
The linear automorphisms in $\Gamma$ are arbitrary permutations
of the coordinates and sign-change automorphisms which allow changing the 
signs of two coordinates.
By applying linear automorphisms, we can assume that $2<x\le y\le z$.

Let $\mu:=2\sqrt{\kappa(u)-2} > 0 $. By Lemma~\ref{lem:MainLemma} the
quadratic reflection $Q_z\in\Gamma$ given by
\begin{equation*}
u= \bmatrix x \\ y \\ z \endbmatrix 
\longmapsto \tilde u = 
\bmatrix x \\ y \\ z' \endbmatrix
\end{equation*}
reduces $z$ by more than $\mu$. If $z'\le 2$ then 
\begin{equation*}
\tilde u\sim (-x,-y,z')\in (-\infty,2]^3,  
\end{equation*}
completing the proof. Otherwise, all three coordinates of $\tilde u$
are greater than $2$ so we repeat the process.
Since each repetition decreases $x+y+z$ by more than $\mu$,
the procedure ends after at most $(x+y+z-6)/\mu$ steps.
\end{proof}

\small

\section*{Appendix: Elements of the modular group}
\let\subsection\rk
We describe in more detail the automorphisms $\Gamma$ and interpret
them geometrically in terms of mapping classes of $M$.


Horowitz \cite{Horowitz} determined the group $\Ak$ of polynomial mappings 
$\C^3\longrightarrow\C^3$ preserving $\kappa$.
We have already observed that the {\em linear automorphisms\/} form
the semidirect product $\Sigma \rtimes \Ss_3$ of the group 
$\Sigma\cong\Z/2\times\Z/2$ of sign-changes (see \S\ref{sec:signchanges}) 
and the symmetric group $\Ss_3$ consisting of permutations of the 
coordinates $x,y,z$. Horowitz proved that the automorphism group of 
$(\C^3,\kappa)$ is generated by the linear automorphism group
$$\Ak\cap \o{GL}(3,\C) = \Sigma \rtimes \Ss_3$$
and the {\em quadratic reflection:\/}
\begin{equation*}
\bmatrix x \\ y \\ z \endbmatrix
\longrightarrow \bmatrix yz - x \\ y \\ z \endbmatrix.
\end{equation*}
This group is commensurable with $\Ou$.

We denote this group by $\Gamma$; it is isomorphic to a semidirect
product 
\begin{equation*}
\Gamma \cong \PGLtz \ltimes (\Z/2\oplus\Z/2)
\end{equation*}
where $\PGLtz$ is the quotient of $\GLtz$ by the elliptic involution
(see below) and $(\Z/2\oplus\Z/2)$ is the group $\Sigma$ of sign-changes.

\subsection{A-1 The elliptic involution}
However, $\Ou$ does not act effectively on $\C^3$. 
To describe elements of $\Ou$, we use the isomorphism 
$h\co \Ou\longrightarrow\GLtz$ discussed in \S\ref{sec:structure}.
The kernel of the homomorphism $\Ou \longrightarrow \Ak$
is generated by $h^{-1}(-I)$.

The elliptic involution is a nontrivial mapping class which
acts trivially on the character variety. This phenomenon is
due to the {\em hyperellipticity \/} of the one-holed torus
(as in \cite{Erg}, \S 10.2). The automorphism $\varepsilon$ of 
$\pi$ given by:
\begin{align*}\label{eq:ellinv}
X & \longmapsto  Y X^{-1} Y^{-1} \sim X^{-1} \\ 
Y & \longmapsto  (Y X) Y^{-1} (X^{-1}Y^{-1}) \sim Y^{-1} \notag\\
XY & \longmapsto X^{-1}Y^{-1} \sim (XY)^{-1} \notag
\end{align*}
preserves $K=[X,Y]$ and 
acts on the homology by the element
\begin{equation*}
h(\varepsilon_*) = -I = \bmatrix -1 & 0 \\ 0 & -1 \endbmatrix \in \GLtz 
\end{equation*}
generating the center of $\GLtz$. Furthermore $\varepsilon$ 
acts identically on the characters $(x,y,z)$.
Thus the homomorphism
\begin{equation*}
\GLtz \overset{h^{-1}}\longrightarrow \Ou
\longrightarrow \Ak
\end{equation*}
factors through $\PGLtz := \GLtz/\{\pm\Id\}$.

Note, however, that $\varepsilon^2 = \iota_{K^{-1}}$.
The automorphism 
$$
\iota_{YX}\circ\varepsilon\co \pi\longrightarrow\pi \qquad
X \longmapsto  X^{-1} \qquad 
Y \longmapsto  Y^{-1} \qquad
XY \longmapsto X^{-1}Y^{-1} \sim (XY)^{-1} 
$$
has order two in $\Aut(\pi)$ but does not preserve $K$.

\subsection{A-2 The symmetric group}
Next we describe the automorphisms of $\pi$ which correspond to
permutations of the three trace coordinates $x,y,z$.
Permuting the two generators $X,Y$ gives:
$$X  \longmapsto  Y \qquad
Y \longmapsto  X \qquad
XY \longmapsto  YX \sim XY .$$
This automorphism $P_{(12)}$ sends $K\longrightarrow K^{-1}$. It acts
on characters by:
\begin{equation*}
(P_{(12)})_*\co
\bmatrix x \\ y \\ z \endbmatrix \longmapsto 
\bmatrix y \\ x \\ z \endbmatrix 
\end{equation*}
and on the homology by
$\bmatrix 0 & 1 \\ 1 & 0 \endbmatrix$. 
Another transposition of the character space is defined by the
involution
\begin{align*}
P_{(13)}\co
X & \longmapsto  Y^{-1}X^{-1} \\
Y & \longmapsto  XYX^{-1} \sim Y \\
XY &\longmapsto  X^{-1} 
\end{align*}
which maps $K\longmapsto K^{-1}$. It acts
on characters by:
\begin{equation*}
(P_{(13)})_*\co
\bmatrix x \\ y \\ z \endbmatrix \longmapsto 
\bmatrix z \\ y \\ x \endbmatrix 
\end{equation*}
and on the homology by
$\bmatrix -1 & 0 \\ -1 & 1 \endbmatrix$. 
The involution
\begin{align*}
P_{(23)}\co
X & \longmapsto  Y^{-1}XY \\
Y & \longmapsto  Y^{-1}X^{-1}  \\
XY &\longmapsto  Y^{-1}
\end{align*}
maps $K\longmapsto K^{-1}$ and acts 
on characters by:
\begin{equation*}
(P_{(23)})_*\co
\bmatrix x \\ y \\ z \endbmatrix \longmapsto 
\bmatrix x \\ z \\ y \endbmatrix 
\end{equation*}
and on the homology by
$\bmatrix 1 & -1 \\ 0 & -1 \endbmatrix$.
The composition $P_{(13)}\circ P_{(12)}$ will be denoted $P_{(123)}$
since the composition of the transposition $(12)$ with the
transposition $(13)$ in $\Ss_3$ equals the 3--cycle $(123)$. 
Applying this composition, we obtain
a 3--cycle of automorphisms
\begin{alignat*}{2}
X & \longmapsto  XYX^{-1}&& \longmapsto  XYX^{-1}Y^{-1}(XY)^{-1} \\
Y & \longmapsto  Y^{-1}X^{-1}  && \longmapsto  (XY)X(XY)^{-1} \\
XY & \longmapsto XYX^{-1}Y^{-1}X^{-1}  && \longmapsto 
(XY)(X^{-1}Y^{-1}X)(XY)^{-1}
\end{alignat*}
which preserve $K$, although $P_{(123)}^3$ equals the elliptic
involution $\varepsilon$, not the identity.
The action $P_{(123)}$ on characters is:
\begin{equation*}
\bmatrix x \\ y \\ z \endbmatrix \overset{P_{(123)}}\longmapsto 
\bmatrix z \\ x \\ y \endbmatrix \overset{P_{(123)}}\longmapsto 
\bmatrix y \\ z \\ x \endbmatrix 
\end{equation*}
which is the {\em inverse \/} of the permutation of coordinates given by 
$(123)$.
The action on homology is given by the respective matrices:
\begin{equation*}
h(P_{(123)}) = \bmatrix 0 & -1 \\ 1 & -1 \endbmatrix,\quad
h(P_{(132)}) = \bmatrix -1 & 1 \\ -1 & 0 \endbmatrix.
\end{equation*}
In summary, we have the following correspondence between permutations of
the coordinates and elements of $\PGLtz$ given by $\tau\mapsto h(P_\tau)$:
$$(12)  \longmapsto \pm \bmatrix 0 & 1 \\ 1 & 0 \endbmatrix \qquad
(13)  \longmapsto \pm \bmatrix 1 & 0 \\ 1 & -1 \endbmatrix \qquad
(23)  \longmapsto \pm \bmatrix -1 & 1 \\ 0 & 1 \endbmatrix $$
$$(132)  \longmapsto \pm \bmatrix 1 & -1 \\ 1 & 0 \endbmatrix \qquad
(123)  \longmapsto \pm \bmatrix 0 & 1 \\ -1 & 1 \endbmatrix$$
The isomorphism $\o{GL}(2,\Z/2)\cong\Ss_3$ relates to the differential
(at the origin) of the mappings in the image of
$\Ou\longrightarrow\Ak$ as follows. The origin is the only isolated
point in the $\Ou$--invariant set $\kappa^{-1}(-2)\cap\R^3$, and is
thus fixed by all of $\Ou$.  (The origin corresponds to the quaternion
representation \eqref{eq:quaternionrep}; see \S\ref{sec:realpoints}.)
Therefore taking the differential at the origin gives a representation
\begin{align*}
\Ou & \longrightarrow \GL{3,\C} \\
\phi & \longmapsto d_{(0,0,0)} \phi_*
\end{align*}
whose image is $\Ss_3$. In particular it identifies with the composition
\begin{equation*}
\Ou \overset{h}\longrightarrow \GLtz \longrightarrow 
\o{GL}(2,\Z/2)\cong\Ss_3
\end{equation*}
where the last arrow denotes reduction modulo 2. 
These facts can be checked by direct computation.

\subsection{A-3 A quadratic reflection}
Here is a mapping class corresponding to a reflection preserving
$(x,y)$. The automorphism $Q_z$ of $\pi$ given by:
\begin{align*}
X & \longmapsto  (XYX^{-1}) X       (XY^{-1}X^{-1})  \sim X \\
Y & \longmapsto  (XYX^{-1}) Y^{-1}  (XY^{-1}X^{-1})  \sim Y^{-1} \\
XY &\longmapsto  (XYX^{-1}) XY^{-1} (XY^{-1}X^{-1})  \sim XY^{-1}
\end{align*}
maps $K\mapsto K^{-1}$ and
\begin{equation*}
(Q_z)_*\co
\bmatrix x \\ y \\ z \endbmatrix  \longmapsto
\bmatrix x \\ y \\ xy -z \endbmatrix 
\end{equation*}
\begin{equation*}
h(Q_z) =\bmatrix 1 & 0 \\ 0 & -1 \endbmatrix.
\end{equation*}
Similarly the other quadratic reflections are:
\begin{equation*}
Q_y\co \begin{cases}
X & \longmapsto  (XY) X       (Y^{-1}X^{-1})  \sim X \\
Y & \longmapsto  (XY) X^{-1}Y^{-1}X^{-1}  (Y^{-1}X^{-1})  
\sim X^{-1}Y^{-1}X^{-1} \\
XY &\longmapsto  Y^{-1} X^{-1} = (XY)^{-1}  
\end{cases}
\end{equation*}
which induces
\begin{equation*}
(Q_y)_*\co
\bmatrix x \\ y \\ z \endbmatrix  \longmapsto
\bmatrix x \\ xz - y \\ z \endbmatrix , \quad
h(Q_y)=  \bmatrix 1 & -2 \\ 0 & -1 \endbmatrix  
\end{equation*}
and
\begin{equation*}
Q_x\co \begin{cases}
X & \longmapsto  (XY) Y^2 X       (Y^{-1}X^{-1})  \sim Y^2 X \\
Y & \longmapsto  (XYX^{-1}) Y^{-1} (XY^{-1}X^{-1})  \sim Y^{-1} \\
XY &\longmapsto  (XY^2) (XY)  (Y^{-2}X^{-1}) \sim XY,
\end{cases}
\end{equation*} 
which induces
\begin{equation*}
(Q_x)_*\co
\bmatrix x \\ y \\ z \endbmatrix  \longmapsto
\bmatrix yz - x \\  y \\ z \endbmatrix ,
h(Q_x)=  \bmatrix 1 & 0 \\ 2 & -1 \endbmatrix.  
\end{equation*}
The square of each of these three reflections is the identity element
of $\Aut(\pi)$.  These reflections correspond to the generators
\eqref{eq:gens2congsubgp} of the level--2 congruence subgroup $\GLtz_{(2)}$.

\subsection{A-4 Another involution}
The automorphism $\nu$
\begin{align*}
X & \longmapsto  Y^{-1} \\
Y & \longmapsto  Y X Y^{-1} \sim X \\
XY & \longmapsto   X Y^{-1}  
\end{align*}
preserves $K$, satisfies $\nu^2 = \varepsilon$ and acts by:
\begin{equation*}
\nu_*\co \bmatrix x \\ y \\ z \endbmatrix \longmapsto
\bmatrix y \\ x \\ xy - z \endbmatrix,
\quad h(\nu_*) = \bmatrix 0 & 1 \\ -1 & 0 \endbmatrix,
\end{equation*}
the composition of the transposition $P_{(12)}$ and the quadratic
reflection $Q_z$. Note that $(P_{(12)})_*$ and $(Q_z)_*$ commute
in $\Ak$.

\subsection{A-5 A Dehn twist}
The automorphism $\tau$
\begin{align*}
X & \longmapsto  X Y \\
Y & \longmapsto  Y  \\
XY & \longmapsto   X Y^2  
\end{align*}
preserves $K$, and acts by:
\begin{equation*}
\tau_*\co \bmatrix x \\ y \\ z \endbmatrix \longmapsto
\bmatrix xy - z \\ y \\ x \endbmatrix,
\quad h(\tau) = \bmatrix 1 & 0 \\ 1 & 1 \endbmatrix,
\end{equation*}
the composition  $P_{(13)}\circ Q_z$.

Similarly the Dehn twist $\tau_X$ around $X$
\begin{align*}
X & \longmapsto  X \notag \\
Y & \longmapsto  Y X \notag \\
XY & \longmapsto   X Y X
\end{align*}
preserves $K$, and acts by:
\begin{equation*} 
\bmatrix x \\ y \\ z \endbmatrix \longmapsto
\bmatrix x \\ xy - z \\ y \endbmatrix,
\quad h(\tau) = \bmatrix 1 & 1 \\ 0 & 1 \endbmatrix,
\end{equation*}
the composition  $P_{(23)}\circ Q_z$.

\end{document}\bye